\documentclass
{article}

\usepackage{amsmath}

\usepackage{amsthm}
\usepackage{amsfonts}
\usepackage{amssymb}

\usepackage{graphicx} 
\usepackage{booktabs}
\usepackage{multirow}
\usepackage{geometry} 
\usepackage{tabularx} 
\usepackage{subcaption}
\usepackage{placeins}
\usepackage[table]{xcolor}
\usepackage{listings}
\usepackage{xcolor}  

\usepackage{tikz}

\usepackage{float}

\usepackage{hyperref}

\usepackage[colorinlistoftodos]{todonotes}
\usepackage{array}
\newcolumntype{L}{>{\raggedright\arraybackslash}X} 

\usepackage[style=authoryear,backend=biber]{biblatex}
\usepackage{algorithm}
\usepackage{algpseudocode}
\usepackage{enumitem}

\newtheorem{remark}{Remark}
\newtheorem{theorem}{Theorem}

\newtheorem{proposition}{Proposition}[section]

\usepackage{tcolorbox}

\usepackage{caption}
\tcbset{
    colback=blue!5!white, 
    colframe=blue!75!black, 
    fonttitle=\bfseries,
    coltitle=black, 
    boxrule=0.5mm, 
    arc=2mm, 
    title=Weakly Dependent Indicator Process
}

\usepackage{parskip}
\usepackage{latexsym}

\usepackage{bm}
\usepackage{bbm}
\newcommand{\one}{\mathbbm{1}}

\usepackage{authblk}

\newcommand{\red}[1]{{\color{black}#1}}

\newcommand{\gr}[1]{{\color{black}#1}}

\begin{document}

\title{Triadic structures in multislice  networks}
\date{}

\author[1]{Kevin Ren}
\author[2]{Tara Trauthwein}
\author[1]{Gesine Reinert}

\affil[1]{Department of Statistics, University of Oxford}
\affil[2]{Institut für Mathematische Stochastik, Universität Münster}

\maketitle

\begin{abstract}
Networks provide a popular representation of complex data. Often, different types of relational measurements are taken on the same subjects. Such data can be represented as a {\it multislice network}, a collection of networks on the same set of nodes, with connections between nodes of different layers.

For the analysis of multislice networks, we take inspiration from the analysis of simple networks, for which small subgraphs (motifs) have proven to be useful; motifs are even seen as building blocks of complex networks. A particular instance of a motif is a triangle, and while triangle counts are well understood for simple network models such as Erd\H{o}s–Rényi  random graphs, with i.i.d.\,distributed edges, even for simple multislice network models little is known about triangle counts. 

Here we address this issue  by extending the analysis of triadic structures to multislice Erdős-Rényi networks. Again taking inspiration from the analysis of sparse  Erd\H{o}s–Rényi  random graphs, we show that the distribution of triangles across multiple layers in a multislice Erdős-Rényi network can be well approximated by an appropriate Poisson distribution. This theoretical result opens the door to statistical goodness of fit tests for multislice networks.
\end{abstract}



\section{Introduction}

Networks have  become an  important tool for describing and analysing complex systems throughout social, biological, physical and mathematical sciences. 
However, often the data complexity is such that multiple measurements are taken on the same set of nodes. A prominent example is that of the Florentine family study in \textcite{padgett_robust_2012}, including marriage relations as well as business relations. This type of data can be represented as a collection of networks on the same set of nodes, with connections between nodes of different networks; such a representation is called a {\it multislice network}. 

Recently, there has been significant research into the study of multislice and, more generally, multilayer networks. The concept of multilayer networks, also called multiplex networks, appears in both engineering, as discussed by \textcite{ChangSeligsonEguchi1996},  and sociology, as detailed by \textcite{WassermanFaust1994};
the
terminology goes back at least to  \textcite{gluckman1955judicial} (Chapter 1, p. 19). Yet, there is a lack 
 of analytic tools for such objects. Although there are many tools available for the statistical analysis of single-layer networks, extending these insights and concepts to multilayer networks remains challenging. Foundational work on multilayer networks by \textcite{Kivela:2014jn}, has provided a framework to study multilayer networks\gr{, and \textcite{Bianconi:2018mn} gives an excellent introduction}; yet there are many gaps in our knowledge regarding the behaviour of multilayer network models.

In particular, counts of small subgraphs (motifs) are among the most important tools for understanding the structural properties of networks, often serving as summary statistics and for comparing different networks see  for example \textcite{alon1997subgraph} and \textcite{alon2007network}. One type of motif, triadic relations, which describe the simplest and most fundamental form of transitivity in a network, are frequently used in social network analysis, see for example \textcite{WassermanFaust1994}.  

As observed in \textcite{picard2008assessing}, in order to assess whether a motif count is exceptional, it is imperative to have a suitable null distribution to compare against. \gr{Motif counts have also been used to understand higher order structures; see for example \textcite{lotito2022higher}. Empirical analysis of multislice networks based on clustering coefficients has been carried out for example in \textcite{bartesaghi2022tensor} and \textcite{bartesaghi2023clustering}. Clustering coefficients are based on counts of triangles and of 2-stars. Yet, a rigorous probabilistic analysis even of triangle counts in simple multislice networks is not available. This paper addresses this gap; it}
focuses on the distribution of triangles within a  multislice Erd\H{o}s–Rényi network (MSER). 

As triangles do not occur independently even in a simple Erd\H{o}s–Rényi network, with independent and identically distributed edge indicators, there is no easy closed form available for their distribution. For a sparse Erd\H{o}s–Rényi graph $G(n, p)$, 
the distribution of
the number of triangles is well approximated by an appropriate Poisson distribution, see for example \textcite{Barbour:1992pa}.  Here we generalise these results to multislice
Erd\H{o}s–Rényi networks.
To illustrate the complications arising, triangles can now span more than one layer; we distinguish 1D triangles for which all edges are in the same layer, 2D triangles in which edges occur across 2 layers, and 3D triangles for which edges occur across 3 layers.
\gr{Using an example from \textcite{Cozzo:2015sr}, in social networks  nodes could represent individuals, and layers could represent multiple modes of
communication. A 1D triangle reflects communication between 3 individuals A, B, and C, using the same mode (perhaps email), while a 2D triangle occurs for example when A and B communicate via email, whereas separately both A and B communicate with C via Slack. If A communicates with B via email, B and C communicate via Slack, while A and C communicate by having lunch together, there are three different modes of communication involved, which is what a 3D triangle counts.}

\textcite{Cozzo:2015sr} offer a \gr{related} definition of multidimensional triangles through paths and give the mean and variance for the number of triangles but does not derive distributional approximations. 

We illustrate the use of the model as well as the distributional approximation by considering the bi-layer Florentine families networks from \textcite{padgett_robust_2012} as well as Lazega's lawyer networks from \textcite{lazega2001collegial}. We find that the MSER model with the same edge probabilities for both layers cannot be rejected for the Florentine data, whereas the MSER for the lawyer multislice network is rejected, using a test at level 5\%. The Poisson approximation here is more of theoretical interest; in the two examples the bounds obtained are not informative. 

This paper is structured as follows. Section~\ref{sec:background} details the background and notation for triangles in a multislice network as well as their uses. Section~\ref{sec:model} introduces the MSER model (short for multislice Erd\H{o}s–Rényi model). Section~\ref{subsec:main} states the  multivariate Poisson approximation, with bounds in total variation distance, and gives an outline of its proof; the proof itself is deferred to Appendix~\ref{app:proofs}.  The use of the result is illustrated in Section~\ref{sec:MC}. The paper ends with a conclusion in Section~\ref{sec:conc}. Appendix~\ref{app:stein} gives more details on Stein's method for multivariate Poisson distributions; Appendix~\ref{app:proofs} contains the detailed proofs of the results in the main text. Python code for the triangle counts is available at \url{https://github.com/rentk/TriCounts}.

\section{Background}\label{sec:background}
\begin{figure}
    \centering
    \includegraphics[width=0.6\linewidth]{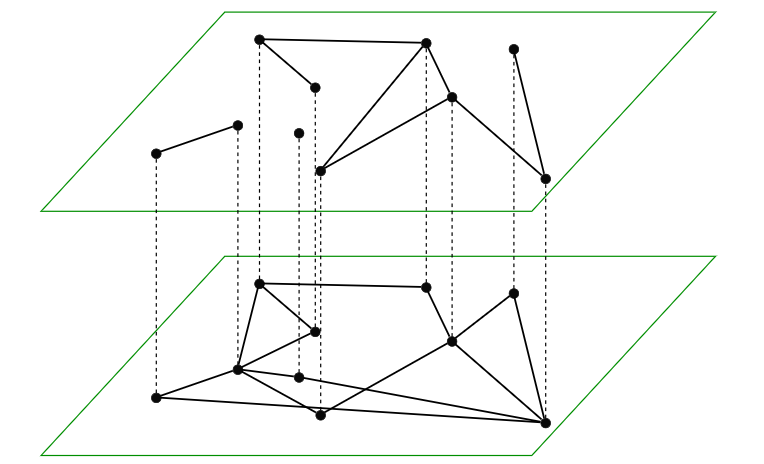}
    \caption{Example of multislice network with two layers. Each layer consists of a Erd\H{o}s–Rényi random graph with 10 nodes and edge probabilities of $0.2$ and $0.3$ respectively. Here  all nodes are connected between layers.}
    \label{fig:mslice}
\end{figure}

A \textit{multislice network} with $L$ slices is  a set of graphs $\{G^i , i = 1, ..., L\}$;  a graph $G^i= (V^i, E^i)$ in layer $i$ has node set 
$V^i = \{u^i | u \in \{1, 2, ..., n\}\}$,
and $E^i$ denotes the edge set in layer $i$. We use the setting of multislice networks with interlinks as  in \textcite{Bianconi:2018mn} and \textcite{Kivela:2014jn}, where copies $u^i$ of the same node $u$ are present in all layers of the network. We call the set $V = \{1, \dots, n\}$ the set of \textit{basis nodes} {(often just called \textit{nodes} in the following)} and we write $u^i \sim v^j$ if node $u$ in layer $i$ is connected to node $v$ in layer $j$. An example of a multislice network can be seen in Figure~\ref{fig:mslice}.

A multislice network can be described by its supra-adjacency matrix, given by 
\begin{align}
    \hat{\mathcal{A}}=\left(\begin{array}{c|c|c|c}
\mathbf{a}^{[1,1]} & \mathbf{a}^{[1,2]} & \cdots & \mathbf{a}^{[1,L]} \\
\hline 
\mathbf{a}^{[2,1]} & \mathbf{a}^{[2,2]} & \cdots & \mathbf{a}^{[2,L]} \\
\hline 
\vdots & \vdots & \ddots & \vdots \\
\hline 
\mathbf{a}^{[L,1]} & \mathbf{a}^{[L,2]} & \cdots & \mathbf{a}^{[L,L]} 
\end{array}\right)  \label{eq:adj}
\end{align}
where each $\mathbf{a}^{[i, j]}$ is a matrix and $\mathbf{a}_{uv}^{[i,j]} = 1$ if $ u^i \sim v^j$, and $0$ otherwise. We use as shorthand the indicator notation $\mathbf{a}_{uv}^{[i,j]}= \mathbbm{1}(u^i \sim v^j).$

\begin{figure}
    \centering
    \includegraphics[width=0.7\linewidth]{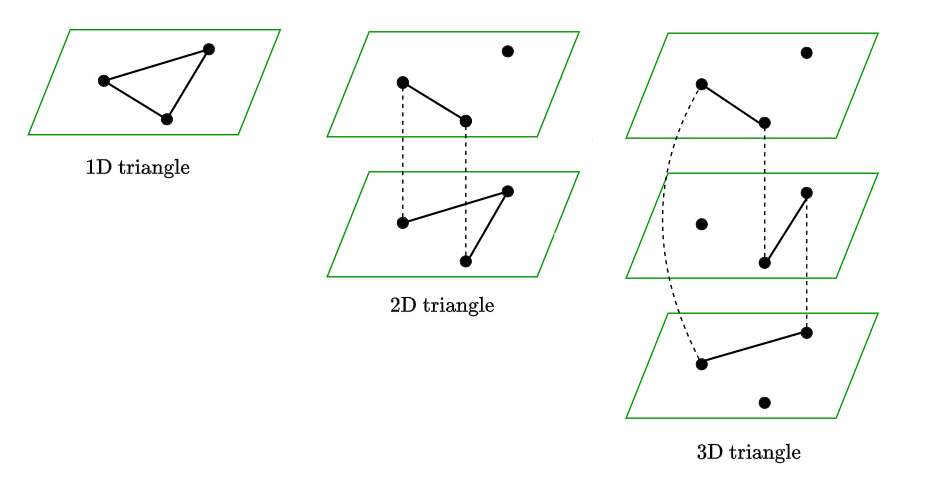}
    \caption{Examples of 1D, 2D and 3D triangles, across one, two or three layers}
    \label{fig:extri}
\end{figure}

Triangles in multislice networks can be of different types, as shown in Figure~\ref{fig:extri}, see also for example \textcite{Cozzo:2015sr}. We call a triangle that is   located in a single layer  a 1D triangle. If the triangle includes edges in two different layers, then we call it a 2D triangle. When edges in three different layers are involved in the triangle, it is called a 3D triangle. Figure~\ref{fig:extri} illustrates the different types of triangles. 

In order to define triangles rigorously, as in \textcite{Kivela:2014jn} we write the adjacency matrix $\hat{\mathcal{A}}$ from \eqref{eq:adj} 
as $ \hat{\mathcal{A}} = {\mathcal{A}} + {\mathcal{C}}$, with  
\begin{align}
\mathcal{A}=\left(\begin{array}{c|c|c|c|c|c}
\mathbf{a}^{[1,1]} & 0 & 0 & \cdots & 0 & 0 \\
\hline 0 & \mathbf{a}^{[2,2]} & 0 & \cdots & 0 & 0 \\
\hline \vdots & \vdots & \vdots & \ddots & \vdots & \vdots \\
\hline 0 & 0 & 0 & \cdots & \mathbf{a}^{[L-1, L-1]} & 0 \\
\hline 0 & 0 & 0 & \cdots & 0 & \mathbf{a}^{[L,L]}
\end{array}\right)
\label{eqn:supraahat}
\end{align}
and
\begin{align}
\mathcal{C}=\left(\begin{array}{c|c|c|c}
\mathbf{0} & \mathbf{a}^{[1,2]} & \cdots & \mathbf{a}^{[1, L]} \\
\hline \mathbf{a}^{[2,1]} & \mathbf{0} & \cdots & \mathbf{a}^{[2, L]} \\
\hline \vdots & \vdots & \ddots & \vdots \\
\hline \mathbf{a}^{[L, 1]} & \mathbf{a}^{[L, 2]} & \cdots & \mathbf{0}
\end{array}\right).
\label{eqn:supraCnew}
\end{align}
The $(nL) \times (nL)$ supra-matrix $\mathcal{A}$  characterizes the intra-layer edges and  the $(nL) \times (nL)$ supra-matrix $\mathcal{C}$ characterizes the inter-layer edges.

In the following, we use $\textbf{Tr}(A)$ to denote the trace of a matrix $A$. 
Using these supra-matrices we can calculate the number of triangles through the use of triadic paths. The number of triadic paths within layer $i$ starting and ending at node $u^i$ is given by $
[\mathcal{A}\mathcal{A}\mathcal{A}]_{u^i,u^i}.
$ 
There are three types of 2D-triadic paths starting and ending at $u^i$ and taking their first step within layer $i$. They are counted by the $u^i,u^i$ entries of $\mathcal{A}\mathcal{A}\mathcal{C}\mathcal{A}\mathcal{C}$, $\mathcal{A}\mathcal{C}\mathcal{A}\mathcal{A}\mathcal{C}$ and $\mathcal{A}\mathcal{C}\mathcal{A}\mathcal{C}\mathcal{A}$. The first of these counts paths that  take two steps within layer $i$, then jump \gr{to another} layer to add an edge in a different layer, and finally jump back to layer $i$, and similarly for the others. The number of 3D-triadic paths is counted in the matrix $\mathcal{A}\mathcal{C}\mathcal{A}\mathcal{C}\mathcal{A}\mathcal{C}$.
These notations characterise triadic paths as walks between intra-layer edges and inter-layer edges.  Thus, in a multislice network, the total numbers $W_1, W_2$ and $W_3$ of 1D, 2D and 3D triangles are 
\begin{align}
    W_1 &= \frac16 {\textbf{Tr}(\mathcal{A}\mathcal{A}\mathcal{A})} \notag\\
    W_2 &= \frac16 \{ {\textbf{Tr}(\mathcal{A}\mathcal{A}\mathcal{C}\mathcal{A}\mathcal{C})}+ {\textbf{Tr}(\mathcal{A}\mathcal{C}\mathcal{A}\mathcal{A}\mathcal{C})} +{\textbf{Tr}(\mathcal{A}\mathcal{C}\mathcal{A}\mathcal{C}\mathcal{A})} \}  \notag\\
    W_3 &= \frac16 {\textbf{Tr}(\mathcal{A}\mathcal{C}\mathcal{A}\mathcal{C}\mathcal{A}\mathcal{C})}, 
\label{eqn:123dtri}
\end{align}
and the total number of triangles is  
\begin{align}
    W = W_1 + W_2 + W_3 .
\end{align}

While the representation \eqref{eqn:123dtri} is useful for computation, in order to disentangle the dependence between triangle counts in a random graph, an alternative representation is useful. 
To this purpose we introduce the notion of \textit{graph isomorphisms}. Given two simple graphs $G=(V(G),E(G))$ and $H=(V(H),E(H))$, an isomorphism of $G$ and $H$ is a bijection $f:V(G) \rightarrow V(H)$ such that $u \sim v \in E(G)$ if and only if $f(u) \sim f(v) \in E(H)$; the graphs $G$ and $H$ are then called isomorphic (see e.g. \textcite[Def. 1.1.20.]{West2001}).

With this definition, exactly $6$ copies of every triadic path are isomorphic (we pick which of the three intra-layer edges we count first, then choose which of its endpoints to start at). \gr{Two paths are thus regarded as equivalent when they differ only by a cyclic change of starting point or by reversal of the direction of traversal, with the node and layer labels fixed.} 

\gr{We}
consider the inter-layer edges as part of the path. 
A \textit{triangle index} is the equivalence class of a corresponding triadic path with respect to isomorphisms. A representative of a triangle index is denoted by $\alpha=(\alpha_1^i,\alpha_2^j,\alpha_3^k)$, by which we mean the triadic path
\begin{align}
\alpha_1^i - \alpha_2^i - \alpha_2^j - \alpha_3^j - \alpha_3^k - \alpha_1^k - \alpha_1^i; 
\label{eq:triapath}
\end{align}
 we remove steps between copies of the same node if we stay in the same layer.
 
 Note that for every equivalence class of triadic paths, we can fix a unique representative. Indeed, both 1D and 3D triangles are uniquely fixed once we go through the nodes in increasing order, and a 2D triangle is unique if we start with the two nodes defining the single edge, setting $\alpha_1$ to be the one with \gr{smaller label
}.
  With this in mind, we can now introduce the index sets $\Gamma_1, \Gamma_2, \Gamma_3$ for all possible 1D, 2D and 3D triangles respectively, as follows.
 \begin{align*}
    &\Gamma_1 := \{\alpha = \overline{(\alpha_1^i,\alpha_2^i,\alpha_3^i)}: \alpha_1<\alpha_2<\alpha_3;\ i \in \{1,\hdots,L\}\}\\
    &\Gamma_2 := \{\alpha = \overline{(\alpha_1^i,\alpha_2^j,\alpha_3^j)}: \alpha_1<\alpha_2;\ \alpha_3 \not\in \{ \alpha_1,\alpha_2\};\ i,j\in \{1,\hdots,L\},\ i\neq j\}\\
    &\Gamma_3:=\{\alpha = \overline{(\alpha_1^i,\alpha_2^j,\alpha_3^k)}: \alpha_1<\alpha_2<\alpha_3;\ i,j,k\in \{1,\hdots,L\},\ i\neq j,\ k \neq i,j\}.
\end{align*}
Here we use the notation $\overline{(\alpha_1^i,\alpha_2^j,\alpha_3^k)}$ to denote the equivalence class of triadic paths defined as in display \eqref{eq:triapath}. We will often abuse notation and simply write $\alpha = (\alpha_1^i,\alpha_2^j,\alpha_3^k)$; moreover, sometimes we refer to $\alpha$ as a triangle.

Given the choice of three nodes $\alpha_1,\alpha_2,\alpha_3$, we can construct $L$ 1D-triangle indices (pick one among $L$ layers), $\binom{L}{2} \cdot 6$ 2D-triangle indices (pick two layers, then decide which of the two layers contains the within-layer  edge, then assign the 3 nodes to the layers) and $\binom{L}{3}\cdot 6$ 3D-triangle indices (pick three layers, then allocate a unique layer to each edge). Hence we have
\begin{align}
    |\Gamma_1| = \binom{n}{3}  L, \quad |\Gamma_2|  =  6 \binom{n}{3} \binom{L}{2}  \text{ and } |\Gamma_3|= 6  \binom{n}{3}  \binom{L}{3} . 
\label{eqn:expgamsetsize}
\end{align}
For each possible triangle index $\alpha = (\alpha_1^i, \alpha_2^j, \alpha_3^k) \in \Gamma_1 \cup \Gamma_2 \cup \Gamma_3$, we define the indicator $X_\alpha$ of its presence in the graph by
\begin{align*}
X_\alpha = 
\begin{cases}
    \one(\alpha_1^i \sim \alpha_2^i) \one(\alpha_1^i \sim \alpha_3^i) \one(\alpha_2^i \sim \alpha_3^i),
    &\text{if } \alpha \in \Gamma_1 \\
    \one(\alpha_1^i \sim \alpha_2^i) \one(\alpha_2^i \sim \alpha_2^j) \one(\alpha_2^j \sim \alpha_3^j) \one(\alpha_3^j \sim \alpha_1^j) \one(\alpha_1^j \sim \alpha_1^i),
    & \text{if } \alpha \in \Gamma_2\\
    \one(\alpha_1^i \sim \alpha_2^i) \one(\alpha_2^i \sim \alpha_2^j) \one(\alpha_2^j \sim \alpha_3^j) \one(\alpha_3^j \sim \alpha_3^k) \one(\alpha_3^k \sim \alpha_1^k) \one(\alpha_1^k \sim \alpha_1^i),
    & \text{if } \alpha \in \Gamma_3.
\end{cases}
\end{align*}

Thus we have the alternative representation to \eqref{eqn:123dtri}, \begin{align}\label{eqn:Wj}
W_1 = \sum_{\alpha \in \Gamma_1} X_\alpha; \quad W_2 = \sum_{\alpha \in \Gamma_2} X_\alpha; \quad
W_3 = \sum_{\alpha \in \Gamma_3} X_\alpha. \end{align}

\section{A multislice Erd\H{o}s–Rényi  model for multislice networks} \label{sec:model}

\gr{In our model,}
each layer $G^i=(V^i,E^i)$ of the network is an Erd\H{o}s-Rényi graph, having independent edge indicators, with edge probabilities   $p_i$, $i=1, \ldots L$.
In this model, two copies of $u^i$ and $u^j$ of the same node in different layers are connected with probability $q$. Thus the \textit{intra-layer} edge probability between two nodes $u^i,v^i$ in layer $i$ is $p_i$ and the \textit{inter-layer} edge (or \textit{down} edge) probability between the same node $u$ and different layers is $q$. 

In this model, each triangle indicator $X_\alpha$ is a Bernoulli random variable with
\begin{equation}
\mathbb{P}(X_\alpha = 1) =
\begin{cases}
    p_i^3   & \text{if } \alpha = (\alpha_1^i,\alpha_2^i,\alpha_3^i) \in \Gamma_1 \\
    p_ip_j^2q^2 & \text{if } \alpha = (\alpha_1^i,\alpha_2^j,\alpha_3^j) \in \Gamma_2 \\
    p_ip_jp_kq^3 & \text{if } \alpha = (\alpha_1^i,\alpha_2^j,\alpha_3^k) \in \Gamma_3.
\end{cases}\label{eqn:pxalpha}
\end{equation}
From \eqref{eqn:pxalpha} and \eqref{eqn:Wj} it follows that 
\begin{align}
\lambda_1&=\mathbb{E} W_1= \sum_{\alpha \in \Gamma_1} \mathbb{E}[X_\alpha] = \sum_{i = 1}^L \binom{n}{3} p_i^3\notag \\
\lambda_2&=\mathbb{E} W_2= \sum_{\alpha \in \Gamma_2} \mathbb{E}[X_\alpha] = 3 \sum_{i = 1}^L \sum_{\substack{j = 1 \\ j \neq i}}^L \binom{n}{3}  p_i p_j^2 q^2 \notag \\
\lambda_3&=\mathbb{E} W_3= \sum_{\alpha \in \Gamma_3} \mathbb{E}[X_\alpha] =  \sum_{i = 1}^L \sum_{\substack{j = 1 \\ j \neq i}}^L
\sum_{\substack{k = 1 \\ k \neq i, j}}^L
\binom{n}{3}  p_i p_j p_k q^3.
\label{eq:means}
\end{align}
\gr{The triangle counts even within  a single layer are correlated; see for example \cite{Barbour:1992pa}.} 

\gr{W}e can give bounds on the covariances  
\begin{align}
        R_{i,j} = \sum_{\alpha\in \Gamma_i} \sum_{\substack{\beta \in \Gamma_j \\ \beta\neq \alpha}} \operatorname{Cov}(X_\alpha,X_\beta) 
     \label{eq:cov}
\end{align} between triangle counts on triangle indices $\alpha$ and $\beta$ that do not span the same triadic path; \gr{see}  \eqref{eq:triapath}.
Before we start with the proof,  we introduce some notation to simplify the presentation. For $m,L \in \mathbb{N}$, let $[L]:=\{1,2,\hdots,L\}$ and denote by $[L]^{m,\neq}$ the set of ordered tuples $(i_1,i_2,\hdots,i_m)$, where $i_k \in [L]$ and all entries are distinct. In Appendix~\ref{app:proofs} we shall prove the following result. 

\medskip 
\begin{proposition}\label{prop:RijBounds}
    We have the following bounds:
    for covariances involving 1D triangles, 
   \begin{align*}
R_{1,1} & = \binom{n}{3}3(n-3) \sum_{i=1}^L p_i^5 (1-p_i)  \leq \frac{1}{2}n^4 \sum_{i=1}^L p_i^5; \\
R_{2,1} & = \sum_{i = 1}^L \sum_{\substack{j = 1 \\ j \neq i}}^L 3\binom{n}{3}q^2\big((n-2)p_i^3p_j^2 {(1-p_i)} + 2(n-3)p_ip_j^4 (1-p_j) + p_ip_j^3 {(1-p_j^2)}\big)\\
&\leq \frac{1}{2}n^3\sum_{i = 1}^L \sum_{\substack{j = 1 \\ j \neq i}}^L q^2\big((n-2)p_i^3p_j^2 + 2(n-3)p_ip_j^4 + p_ip_j^3\big),
\end{align*} 
and
\[
R_{3,1} \leq \frac{1}{2}n^4 \sum_{(i,j,k) \in [L]^{3,\neq}} p_i^3p_jp_kq^3.
\]
For covariances involving 2D triangles but not 1D triangles,

\begin{align*}
    R_{2,2} & \leq \frac{1}{\red{2}}n^3 \sum_{(i,j)\in[L]^{2,\neq}}(4p_i^2p_j^2q^3 + p_i^3p_j^3q^2(1-q^2))
    + \red{4}n^3 \sum_{(i,j,k)\in[L]^{3,\neq}} p_ip_j^2p_kq^4 \\
    &+ \frac{1}{\red{2}}n^4 \sum_{(i,j)\in[L]^{2,\neq}} (8p_i^3p_j^2q^3 + 2p_i^3p_j^3q^3(1-q) + p_ip_j^4q^2 + 4p_i^2p_j^4q^3(1-q) + p_i^3p_j^3q^2(1-q^2)) \\
    &+ \frac{1}{\red{2}}n^4 \sum_{(i,j,k)\in[L]^{3,\neq}} (5p_ip_j^2p_k^2q^4 + 4p_ip_j^3p_kq^4)
    + \red{n^5} \sum_{(i,j)\in[L]^{2,\neq}} (p_i^2p_j^4q^3(1-q) + p_i^3p_j^3q^3(1-q)); \\
    R_{2,3}
    \leq& n^3 \sum_{(i,j,k)\in [L]^{3,\neq}} (2p_ip_j^2p_k q^4 + p_i^2p_j^3p_kq^4(1-q))
    + \frac{3}{2}n^3 \sum_{(i,j,k,\ell)\in[L]^{4,\neq}} p_ip_jp_kp_\ell^2q^5 \\
    &+ n^4 \sum_{(i,j,k) \in [L]^{3,\neq}} (2p_i^2p_j^2p_kq^\red{4} + 2p_i^2p_j^3p_kq^4(1-q) + p_i^3p_jp_kq^4)\\
    &+ \frac{3}{2}n^4 \sum_{(i,j,k,\ell)\in[L]^{4,\neq}} p_ip_jp_k^2p_\ell q^5
    + \frac{1}{4}n^5 \sum_{(i,j,k)\in[L]^{3,\neq}} 4p_ip_j^2p_k^3q^4(1-q).
\end{align*}
For the covariances of 3D triangles,
\begin{align*}
R_{3,3}
&\leq \frac{1}{2}n^3 \sum_{(i,j,k)\in [L]^{3,\neq}} p_ip_j^2p_k^2q^5
+ \frac{1}{2}n^3 \sum_{(i,j,k,\ell) \in [L]^{4,\neq}} (3p_ip_jp_kp_\ell q^5 + p_i^2p_j^2p_kp_\ell q^5(1-q)) \\
&+ \frac{1}{2}n^4 \sum_{(i,j,k) \in [L]^{3,\neq}} (2p_i^2p_j^2p_k^2 q^5(1-q) + 2p_i^2p_jp_k^2 q^4) \\
&+ \frac{1}{2}n^4 \sum_{(i,j,k,\ell)\in [L]^{4,\neq}} (\red{2}p_i^2p_j^2p_kp_\ell q^5(1-q) + \red{4}p_ip_j^2p_kp_\ell q^5)
+ n^3 \sum_{(i,j,k,\ell,m) \in [L]^{5,\neq}} p_ip_jp_kp_\ell p_m q^6 \\
&+ \frac{1}{2} n^5 \sum_{(i,j,k)\in[L]^{3,\neq}} p_i^2p_j^2p_k^2 q^5(1-q) 
+ n^4 \sum_{(i,j,k,\ell,m) \in [L]^{5,\neq}} p_ip_jp_kp_\ell p_m q^6 \\
&+ \frac{1}{4} n^5 \sum_{(i,j,k,\ell)\in[L]^{4,\neq}} 2p_ip_j^2p_k^2p_\ell q^5(1-q).
\end{align*}
\end{proposition}
\begin{remark}\label{rem:sparse}
In the sparse regime in which all $0 \le p_i= c_i/n \le 1$ for some collection of fixed $c_i$'s and fixed $q$'s, all covariances are of the order (at most) $L^\red{5} \max( c_i, 1)^\red{6}/n$ and thus vanish as $n \rightarrow \infty.$ 
\end{remark}

\medskip 
While in the statement of Proposition~\ref{prop:RijBounds} we mostly only give upper bounds on the covariances, the exact expression for $R_{2,1}$ illustrates that triangle counts in different layers are correlated. However, in the sparse regime, this dependence will be weak. The next section exploits this observation. 


\section{A Poisson approximation} \label{subsec:main}

As Proposition~\ref{prop:RijBounds} shows, triangle counts $X_\alpha$ and $X_\beta$ are not independent of each other; as soon as $\alpha$ and $\beta$ share a potential edge, there is dependence. However the dependence is local, in the sense that triangle indicators which do not share any edge indicator are indeed independent. As long as this local dependence structure is weak enough, the counts are however approximately independent, and in the sparse regime  each type of triangle count follows approximately a Poisson distribution; moreover, these Poisson distributions are independent for the different triangle counts. 

To make this intuition precise, we use the following notation. The Poisson distribution with parameter $\lambda$ is denoted by $\rm{Po} (\lambda)$ so that 
$ \rm{Po} (\lambda) \{ k \} = e^{-\lambda} \lambda^k/(k!)$, for $k=0, 1, \ldots$. Given three independent Poisson variables $Z_1,Z_2,Z_3$ with parameters $\lambda_1,\lambda_2,\lambda_3$ respectively, we denote the distribution of the vector $(Z_1,Z_2,Z_3)$ by $\prod_{j=1}^3 \text{Po}(\lambda_j)$. We write ${\cal{L}} (Y)$ for the distribution (or {\it law}) of a random element $Y$. The {\it total variation distance} $d_{TV}$ between two  distributions $P$ and $Q$ on $\{0, 1, \ldots\}^3$ is 
$$d_{TV} (P,Q) = \sup_{A \subset \{0, 1, \ldots\}^3} | P(A) - Q(A)|.$$
Thus, if  $\mathcal{L} (X_1, X_2, X_3) = P$, $\mathcal{L} (Y_1, Y_2, Y_3) = Q$,  and $d_{TV}(P,Q) < \epsilon$ then for all $(a_1, a_2, a_3)\in \{0, 1, \ldots\}^3 $ we can bound 
$$ \mathbbm{P} (Y_i \le a_i, i=1, 2, 3) - \epsilon \le \mathbbm{P} (X_i \le a_i, i=1, 2, 3) \le \mathbbm{P} (Y_i \le a_i, i=1, 2, 3) + \epsilon.$$
If $Q$ is much easier to compute than $P$, then such a bound is useful for example for quantifying the uncertainty of estimates. 
With this notation we have the following result.

\medskip 
\begin{theorem}\label{theorem1}
 The total variation distance between the  joint distribution of triangles $(W_1, W_2, W_3)$ in a MSER
 network, with intra-layer edge probability $p_i$ in layer $i$ and inter-layer edge probability $q$, and a multivariate Poisson distribution $\prod_{j=1}^3 \text{Po}(\lambda_j)$  with $\lambda_i$ defined in \eqref{eq:means} can be bounded as follows:
\begin{align}\label{eq:mainBound}
d_{T V}\left(\mathcal{L}\left(W_1,W_2,W_3\right),  \prod_{j=1}^3 \text{Po}(\lambda_j)\right) \notag 
& \leq {\sum_{i=1}^L} {n \choose 3} 
p_i^6+ 3{\sum_{i=1}^L \sum_{\substack{j=1 \\ j \neq i}}^L} {n \choose 3} 
p_i^2p_j^4 q^4 \notag\\
&+ {\sum_{i=1}^L \sum_{\substack{j=1 \\ j \neq i}}^L \sum_{\substack{k=1 \\ k \neq i,j}}^L}  {n \choose 3} 
p_{i}^2p_{j}^2p_{k}^2  q^6 \notag \\
&+  R_{1,1} + R_{2,2} + R_{3,3} + 2R_{2,1} +  2 R_{3,1}+2 R_{3,2}
\end{align}
where as in Proposition~\ref{prop:RijBounds}, $R_{i,j}= \sum_{\alpha \in \Gamma_i}\sum_{\substack{\beta \in \Gamma_j \\ \beta \neq \alpha}}\operatorname{Cov}(X_\alpha, X_\beta)$. 

If $p_i = p$ for all $i =1, \ldots, L$ and if $q=1$ then we have a simpler bound, namely  
\begin{align}\label{eq:second}
d_{T V}\left(\mathcal{L}\left(W_1,W_2,W_3\right),  \prod_{j=1}^3 \text{Po}(\lambda_j)\right) \leq  L^5 n^4 p^5 + \red{\frac{2}{3}} L^4 n^3 p^4.
\end{align}
\end{theorem}
We note that if $L$ is fixed then the bound in \eqref{eq:second} tends to 0 \gr{not only} in the sparse regime, when $p=c/n$ for some fixed $c>0$, \gr{but also in the  regime $p = o( n^{-4/5})$ which includes a moderately sparse regime}. A similar comment applies to the bound \eqref{eq:mainBound}. 

\medskip 
\begin{remark}
    \gr{Theorem~\ref{theorem1} is a finite-sample result, so in an application its right-hand side can be evaluated explicitly for given number of nodes, number of layers, and edge probabilities. When the bound is small, then conservative tests can be devised based on the approximating Poisson distribution: under the null hypothesis, the distribution of the test statistic in the observed network should behave
    like the approximating Poisson distribution, plus or minus an error bound which is explicit. 
    As a simple example, for constructing conservative confidence intervals, with $\epsilon$ denoting the bound in \eqref{eq:mainBound}, we have for the number $W_2$ of 2D triangles that 
    $$  \mathbbm{P}  ( a \le Z_2 \le b) - \epsilon \le
    \mathbbm{P} ( a\le W_2 \le b) \le \mathbbm{P}  ( a \le Z_2 \le b) + \epsilon$$
    where $Z_2$ has the easy-to-compute Poisson distribution with mean $\lambda_2$.
    We emphasize that practical usefulness depends on the numerical size of the bound $\epsilon$. When the bound exceeds one, as in the examples below, it is uninformative and Monte Carlo calibration may be  preferable.}
\end{remark}

\medskip 
Here is a brief overview of the proof of Theorem~\ref{theorem1}. 
We  employ Stein's method for multivariate Poisson approximation, as in \textcite{AGG90} and \textcite{Barbour:1992pa},  which provides a highly useful way to find quantitative distributional approximations. Observe first that $k \, \rm{Po} (\lambda) \{k\} = \lambda \, \rm{Po} (\lambda) \{k-1\} $, for all $k=1, 2, \ldots$. 
One can show that 
an integer-valued random variable $W$ is Poisson distributed with parameter $\lambda$ if and only if
\[
\mathbb{E}[\lambda f(W+1)-Wf(W)] = 0
\]
for any function $f$ such that $\mathbb{E}[|Zf(Z)|]<\infty$, where $Z \sim {Po}(\lambda)$. Heuristically, if $\mathbb{E}[\lambda f(W+1)-Wf(W)]$ is close to zero, then $W$ should be close in distribution to the law of $Z$. To formalize this,
one introduces the differential equation
\[
\lambda f(x+1) - xf(x) = h(x) - \mathbb{E}[h(Z)],
\]
with solution $f_h$. Considering as functions $h (x) = \mathbbm{1}( x \in A)$, for some $A \subset \{ 0, 1, \ldots\}$ this allows to rewrite the distance
\[
d_{TV}(W,Z) = \sup_{A \subset \{0,1,\hdots\}} |\mathbb{E}[\one(W \in A) - \one(Z \in A)]| = \sup_{A \subset \{ 0, 1, \ldots\}} |\mathbb{E}[\lambda f_A(W+1) - Wf_A(W)]|,
\]
which often yields convenient bounds. We use a multivariate version of this method, allowing to compare a vector $(W_1,W_2,...,W_r)$ to a multivariate Poisson distribution. In a result by \textcite{Barbour:1992pa}, the bound achieved by Stein's Method is applied to sums of indicator random variables $W=\sum_{\alpha} X_\alpha$ and further refined by partitioning the indicators according to how they compare to each other. In our case, indicator random variables $X_\alpha$ signal the presence of triangle $\alpha$, and evaluating the bound boils down to bounding the covariances between triangles, which is in essence a combinatorial question. Much of the proof of Theorem~\ref{theorem1} is devoted to exploring all cases in which triangles can interact and depend on one another.

\section{Illustration} \label{sec:MC}
To illustrate the use of the model, 
 we conduct a goodness of fit test for the MSER model using two methods: first using a Monte Carlo test, then using Theorem~\ref{theorem1}.

 \subsection{Florentine family data} 
 First we illustrate this test on a standard dataset, that of Florentine marriage and business relationships, from  \textcite{padgett_robust_2012}. Each network has 16 nodes referring to each of the families. In each layer there are 120 possible edges; we find 20 edges in the marriage network and 15 in the business network. Often these two networks are modelled by independent Erd\H{o}s–Rényi  models with different edge probabilities estimated by maximum likelihood. 

If however both networks are generated by an MSER model with edge probability $p$, the maximum likelihood estimate for $p$
is 
 $p = 35/240 = 0.146$; we assume that the inter-layer connection probabilities $q$ equal 1, as there is no reason to assume any other value. 
 
 In the network we find 8 1D triangles in total (3 in marriage, 5 in business) and 15 2D triangles. We run a Monte Carlo test with 999 simulated MSER graphs with 16 nodes and 2 layers. We present the Monte Carlo results below for the number of 1D and 2D triangles as well as for their sum, including the interval between the 2.5\% and the 97.5\% quantiles $q_{0.025}$ and $q_{0.975}$ \gr{to further illustrate the spread},
 and the p-value for the test. 

\begin{table}[h!]
\centering
\begin{tabular}{|c|c|c|c|}
\hline
\textbf{} & \textbf{Florentine Family Counts} & \textbf{$[q_{0.025}, q_{0.975}]$} & \textbf{p-value} \\
\hline
1D & 8 & [0, 9] & 0.059 \\
\hline
2D & 15 & [3, 24] & 0.175 \\
\hline
Sum & 23 & [3,32] & 0.117 \\
\hline
\end{tabular}
\caption{Monte Carlo tests for the Florentine families data}
\label{tab:your_label}
\end{table}

\begin{figure}[H]
    \centering
    \includegraphics[width=0.9\linewidth]{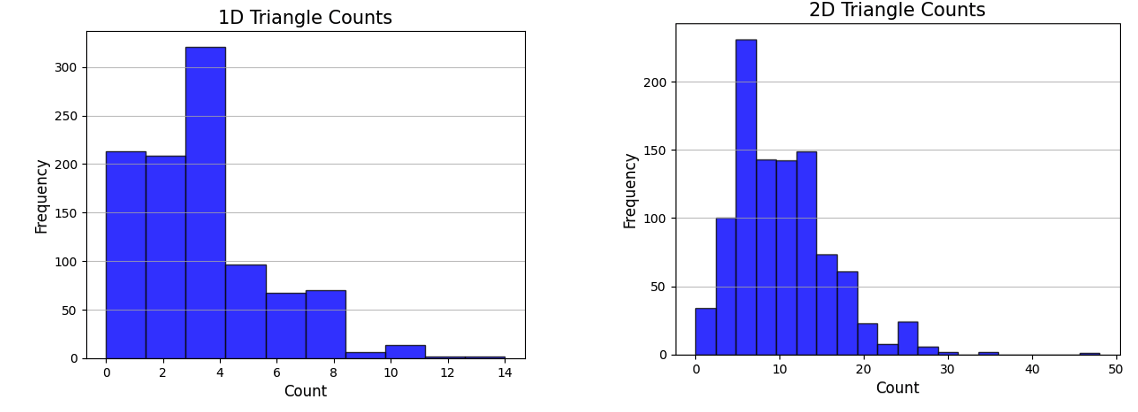}
    \caption{1D and 2D counts from the Monte Carlo simulations, for the Florentine family multislice network}
\end{figure}

In the simulation, we found that there are 37 graphs with \gr{eight}  1D triangles and 40 graphs with at least \gr{nine} 1D triangles. Breaking the 37 ties evenly, the $p$ value was calculated via $(40 + 38/2)/1000  = 0.059$. In our simulations there were 168 graphs with   more than $15$ 2D triangles and  13 graphs with exactly 15 2D triangles.  \gr{As the Florentine family data set is a data set of social interactions, due to a tendency for triadic closure in social networks we would expect to see more triangles than under the MSER model. Hence a one-sided hypothesis test is appropriate. However,}  using a
\gr{one}-sided \gr{test at} 5\% level we do not reject the null hypothesis of the MSER model for any of the three tests. This finding indicates that the two single-layer networks, marriage ties and business ties, may have been generated by a joint mechanism which is reflected in the MSER model.

We can also apply the theoretical bounds  \eqref{eq:second}, with  $n=16$, $p=35/240$,  $q=1$. As for the possible numbers of triangles, we have
$$
|\Gamma| = \binom{16}{3} 2^3  = 4480; \quad |\Gamma_1| = \binom{16}{3}2 = 1120;   \quad |\Gamma_2| = \binom{16}{3} \binom{2}{2} 3! = 3360,
$$
where $\Gamma = \Gamma_1 \cup \Gamma_2$.
The expected total number of triangles is 
$\lambda = |\Gamma|  p^3 = 13.89$. For 1D and 2D triangles, the expected values are $\lambda_1 = 3.47 $ and $ \lambda_2 = 10.42$, respectively; note that as there are only two layers, there are no 3D triangles. Our  bound \eqref{eq:second} on the total variation distance between distribution $\mathcal{L}\left(\left\{W_j\right\}_{j=1}^2\right)$  and a multivariate Poisson distribution $\prod_{j=1}^2 \operatorname{Po}\left(\lambda_j\right)$ with parameters $\lambda_1$ and $\lambda_2$ gives 
$$
d_{T V}\left(\mathcal{L}\left(\left\{W_j\right\}_{j=1}^2\right), \prod_{j=1}^2 \operatorname{Po}\left(\lambda_j\right)\right) 
\leq \red{159}.
$$
This bound is uninformative as the total variation distance $d_{TV}(P, Q)$ always lies within the range $[0,1]$. 
When \red{using the exact expression of the simplified bound derived in the proof of Theorem~\ref{theorem1}}, we can reduce this bound to $\red{47}$; this is still not informative.

\subsection{Lazega's lawyer data}

In \textcite{lazega2001collegial}, relations of different types -- advice, coworker and friendship --  are recorded among 71 lawyers in a New England corporate law firm. Treating the networks as undirected, the advice network has 717 edges, the coworker network has 726 edges, and the friendship network has 399 edges. \gr{The data is available for example at \url{https://networks.skewed.de/net/law_firm}}. The multiplex network  has 5927 1D, 28 440 2D, and 8106 3D triangles. 
Here for an MSER  we take $q=1$ but vary the edge layer probabilities; $p_1= 717/{{71 \choose 2}} ; p_2 = \gr{726/2485} , $ and $p_3=\gr{399/2485}$. 
 Table~\ref{tab:lazega} shows the result from a Monte Carlo test; the null hypothesis of an MSER is clearly rejected as all triadic structure counts are considerably larger than what is seen in simulations from the MSER model.

\begin{table}[h!]
\centering
\begin{tabular}{|c|c|c|c|}
\hline
\textbf{} & \textbf{Lazega's lawyer counts} & \textbf{$[q_{0.025}, q_{0.975}]$} & \textbf{p-value} \\
\hline
1D & 5927 & [2662 , 3430] & 0.001 \\
\hline
2D & 28440 & [13873, 17544] & 0.001 \\
\hline
\gr{3D} & 81\gr{06} & [4096,5269] & 0.001 \\
\hline
\end{tabular}
\caption{Monte Carlo tests for the Lazega lawyer data set}
\label{tab:lazega}
\end{table}

For the multivariate Poisson distribution we obtain the parameters $\lambda_1=  \gr{3034.70} 
, \lambda_2=  \gr{15601.74} 
$, and $\lambda_3= 
\gr{4641.46}
.$ Again the bound on the total variation distance is much larger than 1, making it not informative in this case. We note that this multislice network would not be considered sparse.

\section{Conclusion}\label{sec:conc}

This paper introduces a simple model for multislice networks. As a key network summary, counts of triadic structures are used to assess model fit, both empirically  through Monte Carlo tests as well as theoretically, through a multivariate Poisson approximation. The Poisson approximation is shown to be good in large, sparse networks. For small and dense networks, the bounds are not informative. Yet, they are interesting from a theoretical viewpoint.  

In future work, the distributions of other motif counts could be assessed in a similar fashion. Moreover, the MSER model could be extended to a multiplex stochastic block model setting. For single-layer stochastic block models, approximations for motif counts are available in \textcite{Coulson:2016pa}, see also \textcite{coulson2018compound} for compound Poisson approximations. Similar results should be obtainable for a generalised MSER block model. 

\medskip
{\bf Acknowledgements.}
This research was funded, in part, by UKRI EPSRC grants EP/T018445/1, EP/R018472/1,  EP/X002195/1 and  EP/Y028872/1, as well as by the Research Training Group ``Rigorous Analysis of
Complex Random Systems" (RTG 3027) funded by the Deutsche
Forschungsgemeinschaft (DFG, German Research Foundation – Project
Number 524444762). For the purpose of Open Access, the authors have applied a CC BY public copyright licence to any Author Accepted Manuscript version arising from this submission.

\printbibliography

\appendix

\section{Stein's method for Multivariate Poisson approximation}\label{app:stein}

To show our Theorem~\ref{theorem1}, we use a multivariate Poisson approximation result by Barbour, Holst and Janson (see \textcite{Barbour:1992pa}). Their bound uses the Chen-Stein method and the concept of local dependencies to compare a vector of sums of (possibly dependent) random variables.

To state the result, we need to introduce some notation. Let $\Gamma$ be a set of indices and consider the collection of indicator random variables $(I_\alpha)_{\alpha \in \Gamma}$ with $\mathbb{P}(I_\alpha=1)=:\pi_\alpha$. For any $\alpha \in \Gamma$, we partition the set $\Gamma \setminus \{\alpha\}$ into three subsets $\Gamma_\alpha^-$, $\Gamma_\alpha^+$ and $\Gamma_\alpha^0$, which have the following properties: defining the family of indicators $\left(J_{\beta \alpha}\right)_{\beta \in \Gamma}$ such that 
\[
\mathcal{L}\left( \left(J_{\beta \alpha}\right)_{\beta \in \Gamma} \right) = \mathcal{L} \left( (I_\beta)_{\beta \in \Gamma} | I_\alpha=1 \right),
\]
meaning that the law of $J_{\beta\alpha}$ is the law of $I_\beta$ conditioned on $I_\alpha=1$. The set $\Gamma_\alpha^-$ is such that for every $\beta \in \Gamma_\alpha^-$,
\[
J_{\beta\alpha} \leq I_\beta,
\]
and likewise for $\beta \in \Gamma_\alpha^+$, we have $J_{\beta\alpha} \geq I_\beta$. The $\Gamma_\alpha^0$ contains all other indices.

Assume now in addition that $\Gamma$ can be partitioned into subsets $\Gamma_1,\Gamma_2,\hdots,\Gamma_r$. We are interested in comparing the joint distribution of the sums $W_j = \sum_{\alpha \in \Gamma_j} I_\alpha$ with the distribution of $\Pi_{j=1}^r \text{Po}(\lambda_j)$, which denotes a multivariate Poisson distribution with means $(\lambda_1, \lambda_2, \dots \lambda_r) $, where $\lambda_j:= \mathbb{E} W_j$. The total variation distance between the joint distribution of the $W_j$s and the multivariate Poisson distribution 
$\Pi_{j=1}^r Po(\lambda_j)$ can be bounded as follows.

\medskip 
\begin{theorem}
\label{thm:multpoi}
    [{\textcite[Corollary 10.J.1]{Barbour:1992pa}}] Let $\Gamma=\bigcup_{j=1}^r \Gamma_j$ and $\Gamma\setminus \{\alpha\} = \Gamma_\alpha^+ \cup \Gamma_\alpha^- \cup \Gamma_\alpha^0$, $\left(I_\alpha\right)_{\alpha \in \Gamma}$, as well as $W_j=\sum_{\alpha \in \Gamma_j} I_\alpha$ and $\lambda_j=\mathbb{E} W_j$ be as above. Then
\begin{align}
d_{T V}\left(\mathcal{L}\left(\left\{W_j\right\}_{j=1}^r\right), \Pi_{j=1}^r \text{Po}(\lambda_j)\right) 
&\leq \sum_{\alpha \in \Gamma} \pi_\alpha^2+\sum_{\alpha \in \Gamma} \sum_{\beta \in \Gamma_{\alpha}^{-}}\left|\operatorname{Cov}\left(I_\alpha, I_\beta\right)\right| \notag \\
& +\sum_{\alpha \in \Gamma} \sum_{\beta \in \Gamma_{\alpha}^{+}} \operatorname{Cov}\left(I_\alpha, I_\beta\right)+\sum_{\alpha \in \Gamma} \sum_{\beta \in \Gamma_{\alpha}^0}\left(\mathbb{E} I_\alpha I_\beta+\pi_\alpha \pi_\beta\right).
\end{align}
\label{thm:dtvmpoi}
\end{theorem}

\section{Proofs}
\label{app:proofs}
First we prove Proposition~\ref{prop:RijBounds}. We repeat it here for convenience.

\begingroup
\def\theproposition{\ref{prop:RijBounds}}
\begin{proposition}
For $i,j \in \{1,2,3\}$, recall the notation from Theorem~1:
\[
R_{i,j} = \sum_{\alpha\in \Gamma_i} \sum_{\substack{\beta \in \Gamma_j \\ \beta \neq \alpha}} \operatorname{Cov}(X_\alpha,X_\beta).
\]
We have the following bounds:
For covariances involving 1D triangles, 
\begin{align*}
R_{1,1} & = \binom{n}{3}3(n-3) \sum_{i=1}^L p_i^5 (1-p_i)  \leq \frac{1}{2}n^4 \sum_{i=1}^L p_i^5; \\
R_{2,1} & = \sum_{i = 1}^L \sum_{\substack{j = 1 \\ j \neq i}}^L 3\binom{n}{3}q^2\big((n-2)p_i^3p_j^2 {(1-p_i)} + 2(n-3)p_ip_j^4 (1-p_j) + p_ip_j^3 {(1-p_j^2)}\big)\\
&\leq \frac{1}{2}n^3\sum_{i = 1}^L \sum_{\substack{j = 1 \\ j \neq i}}^L q^2\big((n-2)p_i^3p_j^2 + 2(n-3)p_ip_j^4 + p_ip_j^3\big),
\end{align*} 
and
\[
R_{3,1} \leq \frac{1}{2}n^4 \sum_{(i,j,k) \in [L]^{3,\neq}} p_i^3p_jp_kq^3.
\]
For covariances involving 2D triangles but not 1D triangles,

\begin{align*}
    R_{2,2} & \leq \frac{1}{\red{2}}n^3 \sum_{(i,j)\in[L]^{2,\neq}}(4p_i^2p_j^2q^3 + p_i^3p_j^3q^2(1-q^2))
    + \red{4}n^3 \sum_{(i,j,k)\in[L]^{3,\neq}} p_ip_j^2p_kq^4 \\
    &+ \frac{1}{\red{2}}n^4 \sum_{(i,j)\in[L]^{2,\neq}} (8p_i^3p_j^2q^3 + 2p_i^3p_j^3q^3(1-q) + p_ip_j^4q^2 + 4p_i^2p_j^4q^3(1-q) + p_i^3p_j^3q^2(1-q^2)) \\
    &+ \frac{1}{\red{2}}n^4 \sum_{(i,j,k)\in[L]^{3,\neq}} (5p_ip_j^2p_k^2q^4 + 4p_ip_j^3p_kq^4)
    + \red{n^5} \sum_{(i,j)\in[L]^{2,\neq}} (p_i^2p_j^4q^3(1-q) + p_i^3p_j^3q^3(1-q)); \\
    R_{2,3}
    \leq& n^3 \sum_{(i,j,k)\in [L]^{3,\neq}} (2p_ip_j^2p_k q^4 + p_i^2p_j^3p_kq^4(1-q))
    + \frac{3}{2}n^3 \sum_{(i,j,k,\ell)\in[L]^{4,\neq}} p_ip_jp_kp_\ell^2q^5 \\
    &+ n^4 \sum_{(i,j,k) \in [L]^{3,\neq}} (2p_i^2p_j^2p_kq^\red{4} + 2p_i^2p_j^3p_kq^4(1-q) + p_i^3p_jp_kq^4)\\
    &+ \frac{3}{2}n^4 \sum_{(i,j,k,\ell)\in[L]^{4,\neq}} p_ip_jp_k^2p_\ell q^5
    + \frac{1}{4}n^5 \sum_{(i,j,k)\in[L]^{3,\neq}} 4p_ip_j^2p_k^3q^4(1-q).
\end{align*}
For the covariances of 3D triangles,
\begin{align*}
R_{3,3}
&\leq \frac{1}{2}n^3 \sum_{(i,j,k)\in [L]^{3,\neq}} p_ip_j^2p_k^2q^5
+ \frac{1}{2}n^3 \sum_{(i,j,k,\ell) \in [L]^{4,\neq}} (3p_ip_jp_kp_\ell q^5 + p_i^2p_j^2p_kp_\ell q^5(1-q)) \\
&+ \frac{1}{2}n^4 \sum_{(i,j,k) \in [L]^{3,\neq}} (2p_i^2p_j^2p_k^2 q^5(1-q) + 2p_i^2p_jp_k^2 q^4) \\
&+ \frac{1}{2}n^4 \sum_{(i,j,k,\ell)\in [L]^{4,\neq}} (\red{2}p_i^2p_j^2p_kp_\ell q^5(1-q) + \red{4}p_ip_j^2p_kp_\ell q^5)
+ n^3 \sum_{(i,j,k,\ell,m) \in [L]^{5,\neq}} p_ip_jp_kp_\ell p_m q^6 \\
&+ \frac{1}{2} n^5 \sum_{(i,j,k)\in[L]^{3,\neq}} p_i^2p_j^2p_k^2 q^5(1-q) 
+ n^4 \sum_{(i,j,k,\ell,m) \in [L]^{5,\neq}} p_ip_jp_kp_\ell p_m q^6 \\
&+ \frac{1}{4} n^5 \sum_{(i,j,k,\ell)\in[L]^{4,\neq}} 2p_ip_j^2p_k^2p_\ell q^5(1-q).
\end{align*}
\end{proposition}
\addtocounter{proposition}{-1}
\endgroup

\medskip 
\begin{proof}
We use {the count representation \eqref{eqn:Wj} with}
the sets of indices $\Gamma=\Gamma_1 \cup \Gamma_2 \cup \Gamma_3$ introduced in Section~\ref{sec:background}, and we bound the covariances one by one. 

\subsubsection*{Bound for $R_{1,1}$}

Fix an index $\alpha \in \Gamma_1$ living in layer $i$. Any $X_\beta$ with $\beta$ not sharing a potential edge with $\alpha$ is independent of $X_\alpha$ and the covariance is zero. We also exclude the triangle index $\beta$ which spans the same triadic path as $\alpha$. We thus count all $\beta \in \Gamma_1$ sharing exactly one potential edge with $\alpha$ and having exactly two more potential edges on the same layer $i$. Noting that $ \mathbbm{E} X_\alpha = p_i^3$ this gives 
\[
\operatorname{Cov}(X_\alpha,X_\beta) = 
\mathbbm{E} X_\alpha X_\beta - p_i^6 = 
p_i^5-p_i^6 \leq p_i^5.
\]
There are $\binom{n}{3}$ choices for the three nodes of $\alpha$, and $3(n-3)$ choices for the triangle index $\beta$ once all edges of $\alpha$ are fixed (pick the edge they have in common, then pick an additional node in the same layer to form the other two edges of $\beta$). We thus have:
\begin{align*}
R_{1,1} 
&= \sum_{\alpha \in \Gamma_1} \sum_{\substack{\beta \in \Gamma_1 \\ \beta \neq \alpha}} \sum_{i=1}^L \one(\alpha, \beta \text{ in layer } i \text{ sharing one edge})\cdot(p_i^5 - p_i^6) = \binom{n}{3}3(n-3) \sum_{i=1}^L p_i^5 (1-p_i) \\
&
\leq \binom{n}{3}3(n-3) \sum_{i=1}^L p_i^5.
\end{align*}

\subsubsection*{Bound for $R_{2,1}$}
Fix $\alpha \in \Gamma_2$ and denote by $i$ the layer containing the potential single edge and $j$ the layer containing the other two potential edges. For $R_{2,1}$ we consider covariances with $X_\beta$ for $\beta$ the index of a potential single-layer triangle. Again we count the number of indices $\beta$ sharing a potential edge with $\alpha$. From \eqref{eqn:pxalpha}, $ \mathbbm{E} X_\alpha = p_i p_j^2 q^2.$ There are $n-2$ ways of choosing a 1D triangle index $\beta$ in layer $i$, which then shares a potential edge in layer $i$ with $\alpha$, so that  for such $\beta$ we have  $\operatorname{Cov}(X_\alpha,X_\beta)
{=} p_i^3p_j^2q^2 {-p_i^4p_j^2q^2}$. Similarly, there are $2(n-3)$ ways to get a potential triangle index $\beta$ in layer $j$ which shares one potential edge with $\alpha$, in which case $\operatorname{Cov}(X_\alpha,X_\beta)  
= p_ip_j^4q^2{-p_ip_j^5q^2}$.  Finally, there is  a single way to pick a triangle index $\beta$ which shares two edges with $\alpha$, in which case $\operatorname{Cov}(X_\alpha,X_\beta) =  p_ip_j^3q^2{-p_ip_j^5q^2}$. Thus,
\begin{align*}
    R_{2,1} &= \sum_{i = 1}^L \sum_{\substack{j = 1 \\ j \neq i}}^L 3\binom{n}{3}q^2\big((n-2)p_i^3p_j^2 {(1-p_i)} + 2(n-3)p_ip_j^4 (1-p_j) + p_ip_j^3 {(1-p_j^2)}\big)\\
    &\leq \sum_{i = 1}^L \sum_{\substack{j = 1 \\ j \neq i}}^L 3\binom{n}{3}q^2\big((n-2)p_i^3p_j^2 + 2(n-3)p_ip_j^4 + p_ip_j^3\big).
\end{align*}

\subsubsection*{Bound for $R_{3,1}$}
Each edge in a 3D triangle at $\alpha$ living in layers $i,j,k$ can share an  edge with $(n-2)$ 1D triangles. Say the 1D triangle index $\beta$ lives in layer $i$, then this gives a covariance bound of $Cov(X_\alpha, X_\beta) \leq p_i^3 p_j p_k q^3$. By symmetry we get  
\begin{align*}
R_{3,1} &\leq \binom{n}{3} (n-2) \sum_{i=1}^L\sum_{\substack{j=1\\j\neq i}}^L \sum_{\substack{k=1 \\ k \neq i,j}}^L (p_i^3p_jp_kq^3 + p_ip_j^3p_kq^3 + p_ip_jp_k^3q^3)\\
&= 3\binom{n}{3} (n-2) \sum_{i=1}^L\sum_{\substack{j=1\\j\neq i}}^L \sum_{\substack{k=1 \\ k \neq i,j}}^L p_i^3p_jp_kq^3.
\end{align*}

\subsubsection*{Bound for $R_{2,2}$}
This calculation is considerably more involved;  we split it into different subcases depending on the number of nodes and the number of layers that are added. Fix $\alpha\in \Gamma_2$ and denote the layer with the single edge by $i$ and the one with two edges by $j$. 

\textbf{Case 1:} No new layers, no new nodes\\
When using the same layers and nodes as for $\alpha$, there are six possible 2D triangle indices $\beta$ which share at least one potential edge with $\alpha$, see Figure~\ref{fig:22f0g0}. The black edges indicate a 2D triangle at $\alpha$; the orange edges indicate the potential 2D triangles which share at least one edge with $\alpha$. Although there are no new nodes, there are different ways of distributing the edges to obtain a 2D triangle. We adopt the notation $\{ a, b\}$ where $a$ gives the number of shared intra-layer edges and $b$ the number of shared down edges. For example, in the first panel of the first row of Figure~\ref{fig:22f0g0}, all three intra-layer edges are shared and both down-edges are shared. Hence the triadic paths are identical and thus do not contribute to $R_{2,2}$. {In the second panel, one intra-layer edge and one down edge are shared. For this $\beta$, we get that}
{\[
\operatorname{Cov} (X_\alpha, X_\beta) = p_i^2p_j^3q^3 - p_i p_j^2 q^2  \mathbbm{E} X_\beta \leq p_i^2 p_j^3  q^3 .
\]}
Summing over all possible such $\alpha $ and $\beta$ gives a bound of
\begin{align*}
C_1&:= \red{3}\binom{n}{3} \sum_{(i,j)\in[L]^{2,\neq}} (2p_i^2p_j^3q^3 + 2p_i^2p_j^2q^3 + p_i^3p_j^3q^2(1-q^2)) \\
& \leq \red{3}\binom{n}{3} \sum_{(i,j)\in[L]^{2,\neq}}(4p_i^2p_j^2q^3 + p_i^3p_j^3q^2(1-q^2)).
\end{align*}
\begin{figure}[H]
    \centering
    \includegraphics[width=0.65\linewidth]{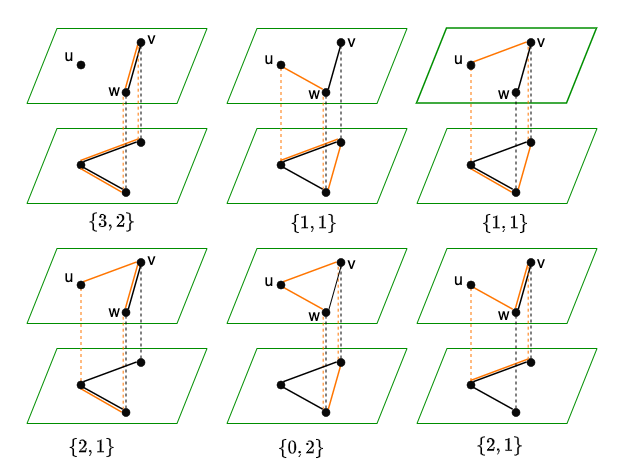}
    \caption{Dependent triangles for Case 1}
    \label{fig:22f0g0}
\end{figure}

\textbf{Case 2:} One new layer, no new nodes\\
We add a new layer $k$, but no new nodes. There are now six possible triangles using layers $j,k$ and three triangles sharing edges in layers $i,k$, see Figures~\ref{fig:22f1g0a} and~\ref{fig:22f1g0b}. Again arguing as in the first case and summing over all possible $\alpha$ and $\beta$ in this case, we obtain 
\begin{align*}
C_2:&= \red{3}\binom{n}{3} \sum_{(i,j,k)\in[L]^{3,\neq}} (2p_ip_j^2p_k^2q^4 + p_ip_j^2 p_kq^4 + 2p_ip_j^3p_kq^4 + p_i^2p_jp_k^2q^4 + 2p_ip_j^2p_k^2q^4)\\
&\leq \red{24}\binom{n}{3} \sum_{(i,j,k)\in[L]^{3,\neq}} p_ip_j^2p_kq^4.
\end{align*}

Note that the third picture in the top row of Figure 5 is a $\{0,0\}$ configuration; in this one,  $\alpha$ and $\beta$ do not share any edges, hence the covariance is zero and this does not contribute to $R_{2,2}$.
\begin{figure}[H]
    \centering
    \includegraphics[width=0.6\linewidth]{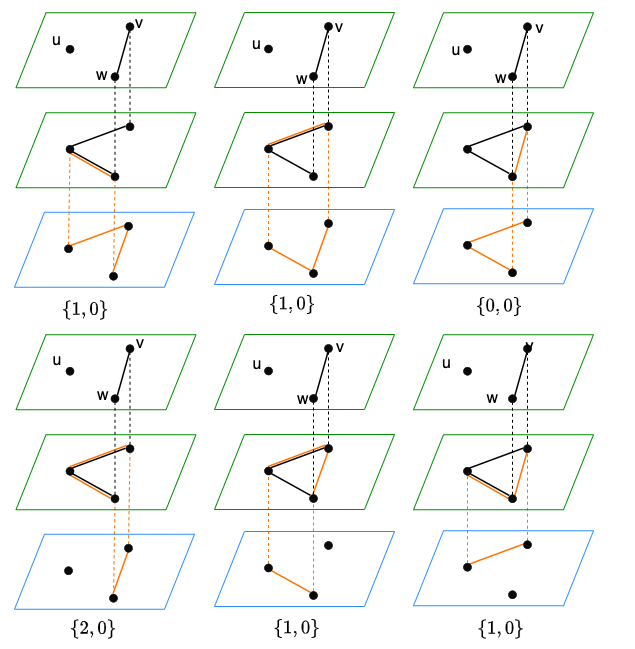}
    \caption{Dependent triangles in Case 2, sharing the layer with two edges}
    \label{fig:22f1g0a}
\end{figure}
\begin{figure}[H]
    \centering
    \includegraphics[width=0.6\linewidth]{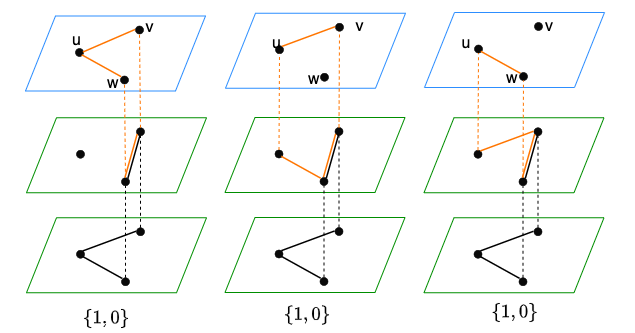}
    \caption{Dependent triangles in Case 2, sharing the layer with one edge}
    \label{fig:22f1g0b}
\end{figure}

\textbf{Case 3:} No new layers, one new node\\
With the addition of a new node, we need to pick two nodes from $\alpha$. There are three choices: picking the node with degree two in layer $j$, and one of the other two nodes ($2$ options), or picking the two nodes adjacent to the single edge. All choices give six possible triangles each illustrated in Figures~\ref{fig:22f0g1a} and \ref{fig:22f0g1b}. 
\begin{figure}[H]
    \centering
    \includegraphics[width=0.6\linewidth]{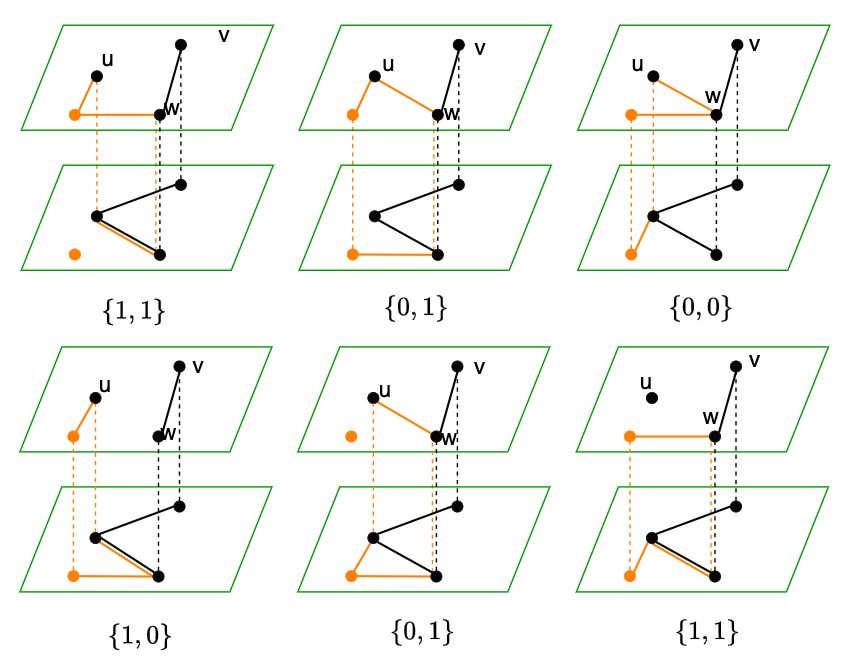}
    \caption{Dependent triangles in Case 3, sharing the isolated node.}
    \label{fig:22f0g1a}
\end{figure}

\begin{figure}[H]
    \centering
    \includegraphics[width=0.6\linewidth]{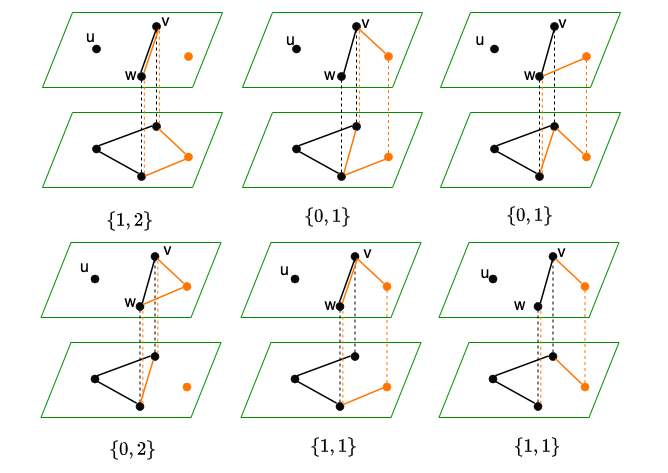}
    \caption{Dependent triangles in Case 3, sharing the isolated edge}
    \label{fig:22f0g1b}
\end{figure}

We get 
\begin{align*}
    C_3:=& \red{3}\binom{n}{3} (n-3) \sum_{(i,j)\in[L]^{2,\neq}} \bigg(2\big(p_i^3p_j^2q^3 + p_i^3p_j^3q^3(1-q) + p_i^2p_j^3q^4 +p_i^2p_j^4q^3(1-q) \\
    &\hspace{3cm}+ p_i^2p_j^3q^3 \big) + p_ip_j^4q^2 + 2p_i^2p_j^4q^3(1-q) + p_i^3p_j^3q^2(1-q^2) + 2p_i^2p_j^3q^3\bigg) \\
    &\leq \red{3}\binom{n}{3} (n-3) \sum_{(i,j)\in[L]^{2,\neq}} \big(8p_i^3p_j^2q^3 + 2p_i^3p_j^3q^3(1-q) + p_ip_j^4q^2 + \\
    &\hspace{8cm}4p_i^2p_j^4q^3(1-q) + p_i^3p_j^3q^2(1-q^2)\big).
\end{align*}

\textbf{Case 4:} One new layer, one new node\\
Here we can only have dependent triangles if an intra-layer edge is shared. The first row in Figure~\ref{fig:22f1g1b} shows the possibilities when sharing one of the two edges in layer $j$ (and there are $2$ choices for the shared edge, so this row counts twice), the second row illustrates sharing the isolated edge. We get 
\begin{align*}
    C_4:=& \red{3}\binom{n}{3}(n-3) \sum_{(i,j,k)\in[L]^{3,\neq}} \bigg(2\big(p_ip_j^2p_k^2q^4 + 2p_ip_j^3p_kq^4\big) + p_ip_j^2p_k^2q^4 + 2p_i^2p_j^2p_kq^4\bigg)\\
    &\leq \red{3}\binom{n}{3}(n-3) \sum_{(i,j,k)\in[L]^{3,\neq}} (5p_ip_j^2p_k^2q^4 + 4p_ip_j^3p_kq^4).
\end{align*}

\begin{figure}[H]
    \centering
    \includegraphics[width=0.6\linewidth]{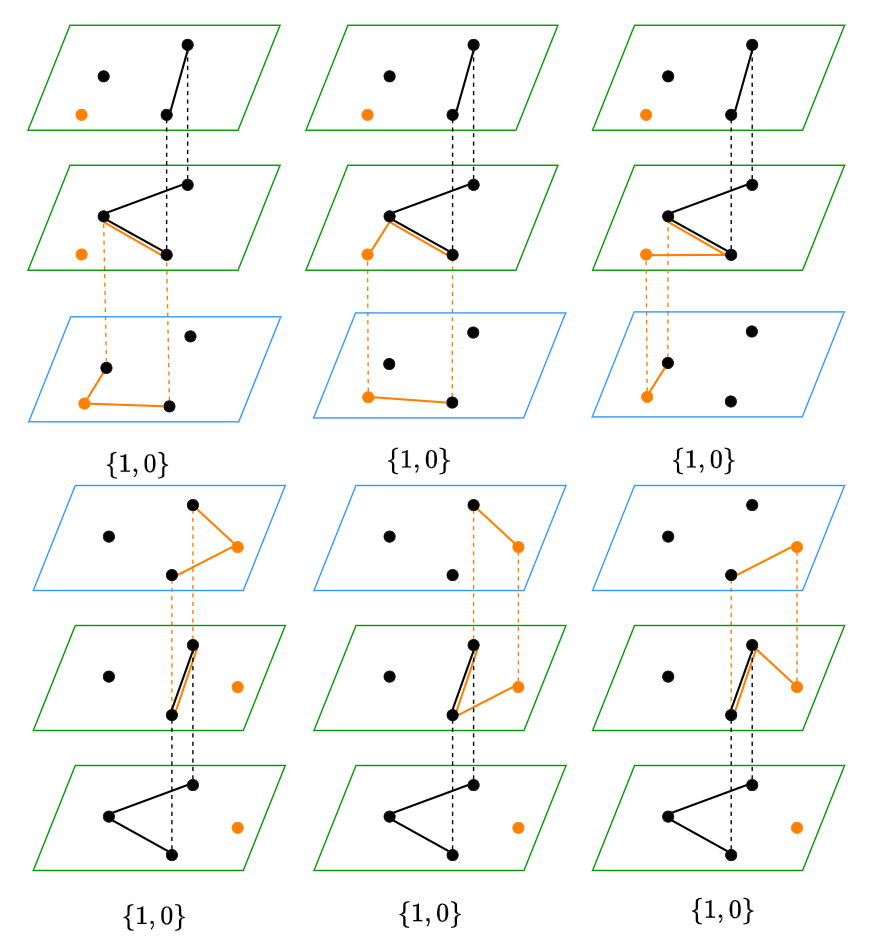}
    \caption{Dependent triangles in Case 4; note that in the second row, the new layer $k$ is displayed in blue at the top.}
    \label{fig:22f1g1b}
\end{figure} 

\textbf{Case 5:} No new layers, two new nodes\\
In this case, only down edges can be shared. There are two choices of down edges, and for each we can construct four triangles $\beta$, see Figure~\ref{fig:22f0g2}. We get 
\[
C_5:= \red{12}\binom{n}{3}\red{\binom{n-3}{2}} \sum_{(i,j)\in[L]^{2,\neq}} (p_i^2p_j^4q^3(1-q) + p_i^3p_j^3q^3(1-q)).
\]
\begin{figure}[H]
    \centering
    \includegraphics[width=0.5\linewidth]{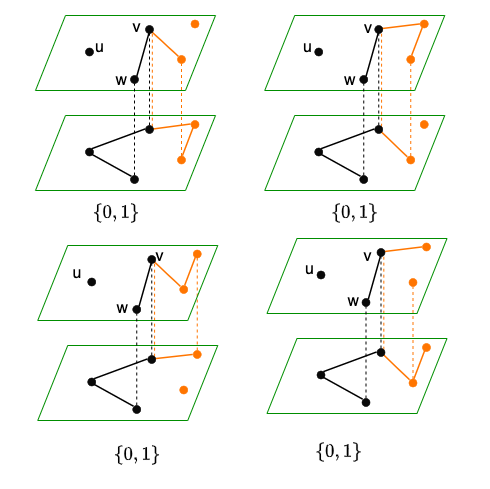}
    \caption{Dependent triangles in Case 5}
    \label{fig:22f0g2}
\end{figure} 

Summing and bounding ${n \choose 3} \le n^3/6$ and $n-3 \le n$ \red{and $\binom{n-3}{2} \leq n^2/2$}  gives that 
\begin{align*}
    R_{2,2} \leq& C_1 + C_2 + C_3 + C_4 + C_5\\
    \leq &\frac{1}{\red{2}}n^3 \sum_{(i,j)\in[L]^{2,\neq}}(4p_i^2p_j^2q^3 + p_i^3p_j^3q^2(1-q^2)) \\
    &+ \red{4}n^3 \sum_{(i,j,k)\in[L]^{3,\neq}} p_ip_j^2p_kq^4 \\
    &+ \frac{1}{\red{2}}n^4 \sum_{(i,j)\in[L]^{2,\neq}} (8p_i^3p_j^2q^3 + 2p_i^3p_j^3q^3(1-q) + p_ip_j^4q^2 + 4p_i^2p_j^4q^3(1-q) + p_i^3p_j^3q^2(1-q^2)) \\
    &+ \frac{1}{\red{2}}n^4 \sum_{(i,j,k)\in[L]^{3,\neq}} (5p_ip_j^2p_k^2q^4 + 4p_ip_j^3p_kq^4) \\
    &+ \red{n^5} \sum_{(i,j)\in[L]^{2,\neq}} (p_i^2p_j^4q^3(1-q) + p_i^3p_j^3q^3(1-q)).
\end{align*}

\subsubsection*{Bound for $R_{3,2}$} 
Here we assume $\alpha \in \Gamma_3$, with edges in layers $i,j,k$, and consider several cases.

\textbf{Case 1:} No new layers, no new nodes\\
To create a 2D triangle $\beta$, we pick two layers out of $i,j,k$, which gives $3$ choices. Now six different triangles can be built, see Figure~\ref{fig:32f0g0}. As before, we define a quantity $B_1$ which will be added to the bound at a later stage. We have:
\begin{align*}
B_1:&= 3\binom{n}{3} \sum_{(i,j,k)\in [L]^{3,\neq}} (2p_ip_j^2p_k q^4 + 2p_i^2p_j^2p_kq^5 + 2p_i^2p_j^3p_kq^4(1-q)) \\
&\leq n^3 \sum_{(i,j,k)\in [L]^{3,\neq}} (2p_ip_j^2p_k q^4 + p_i^2p_j^3p_kq^4(1-q)).
\end{align*}
Note that we used the symmetry of the sum to group terms of the type $p_ip_j^2p_k$ and $p_i^2p_jp_k$ together as $2p_ip_j^2p_k$, and similarly for $p_i^2p_j^3p_k$ and $p_i^3p_j^2p_k$.
\begin{figure}[H]
    \centering
    \includegraphics[width=\linewidth]{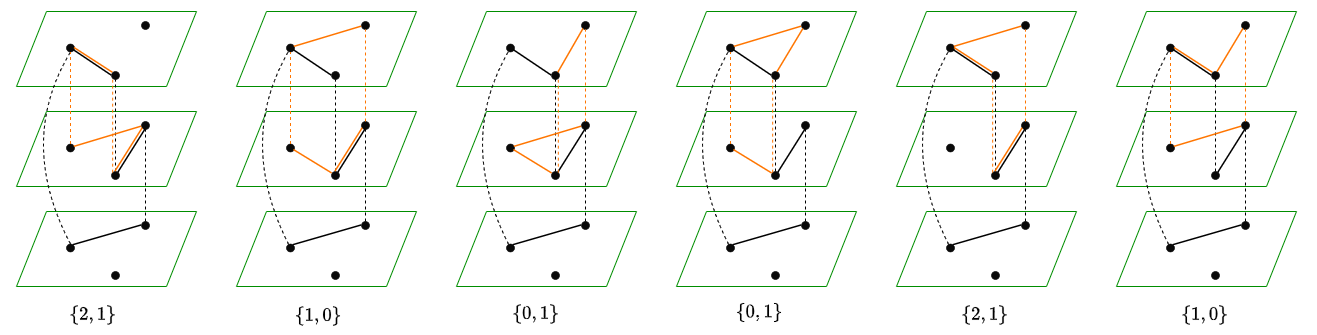}
    \caption{Dependent triangles for Case 1}
    \label{fig:32f0g0}
\end{figure}

\textbf{Case 2:} One new layer, no new nodes\\
Given one new layer, we have $3$ choices to pick one additional layer from $i,j,k$, which gives rise to three possible triangles $\beta$. In Figure~\ref{fig:32f1g0} we show the possible triangles once one layer has been picked. Note that we do not show the other layers containing edges of $\alpha$. We get
\[
B_2:= 3\binom{n}{3} \sum_{(i,j,k,\ell)\in[L]^{4,\neq}} (p_ip_jp_kp_\ell^2q^5 + 2p_ip_jp_k^2p_\ell q^5) \leq \frac{3}{2}n^3 \sum_{(i,j,k,\ell)\in[L]^{4,\neq}} p_ip_jp_kp_\ell^2q^5.
\]
\begin{figure}[H]
    \centering
    \includegraphics[width=0.5\linewidth]{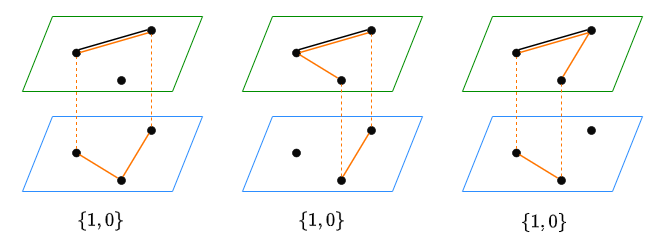}
    \caption{Dependent triangles in Case 2}
    \label{fig:32f1g0}
\end{figure}

\textbf{Case 3:} No new layers, one new node\\
Given a new node, there are $3$ ways in which one can pick the other two nodes, after which there are only $2$ choices of layers which give dependent triangles $\beta$, the third choice resulting in triangles which have no common edges. See Figure~\ref{fig:32f0g1} for an illustration of possible $\beta$ once layers and nodes have been fixed. We have thus:
\begin{align*}
B_3:&= 6\binom{n}{3} (n-3) \sum_{(i,j,k) \in [L]^{3,\neq}} (p_i^2p_j^2p_kq^5 + \red{p_i^2p_j^2p_kq^4 }+ p_i^2p_j^3p_kq^4(1-q) \\
&\hspace{8cm}+ p_i^3p_j^2p_kq^4(1-q) + p_i^3p_jp_kq^4) \\
&\leq n^4 \sum_{(i,j,k) \in [L]^{3,\neq}} (2p_i^2p_j^2p_kq^\red{4} + 2p_i^2p_j^3p_kq^4(1-q) + p_i^3p_jp_kq^4).
\end{align*}

\begin{figure}[H]
    \centering
    \includegraphics[width=\linewidth]{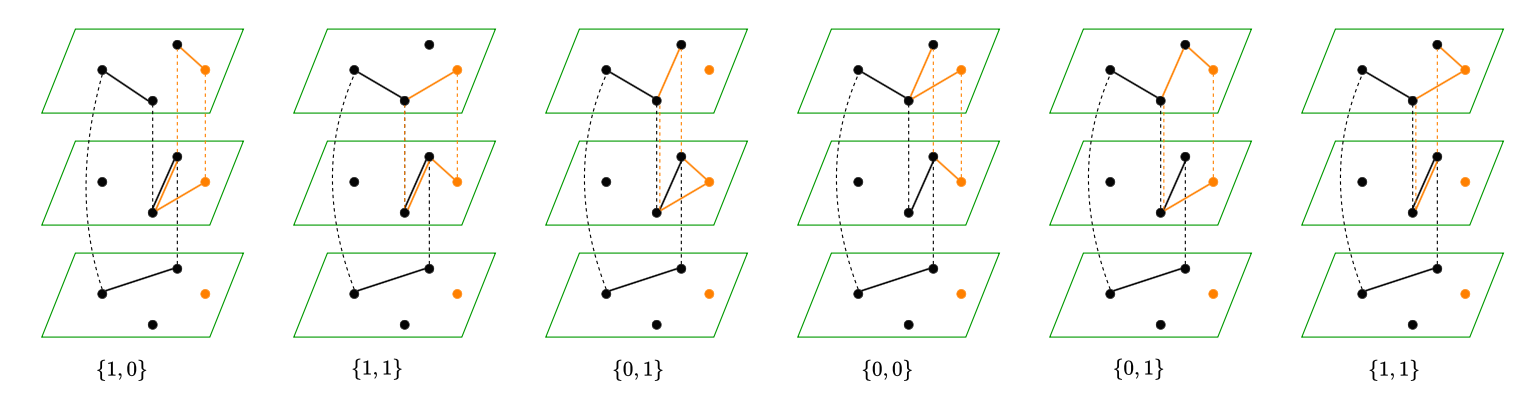}
    \caption{Dependent triangles in Case 3}
    \label{fig:32f0g1}
\end{figure}

\textbf{Case 4:} One new layer, one new node\\
Given a new layer $\ell$, there are $3$ choices of layers among $i,j,k$. Once the layer has been picked, the only way to create a triangle with common edges is to pick the edge of $\alpha$ already present in the layer. Figure~\ref{fig:32f1g1} enumerates the possibilities once a layer has been picked, omitting the other two layers in the picture. We have:
\[
B_4:= 3\binom{n}{3}(n-3) \sum_{(i,j,k,\ell)\in[L]^{4,\neq}} (p_ip_jp_kp_\ell^2q^5 + 2p_ip_jp_k^2p_\ell q^5) \leq \frac{3}{2}n^4 \sum_{(i,j,k,\ell)\in[L]^{4,\neq}} p_ip_jp_k^2p_\ell q^5.
\]
\begin{figure}[H]
    \centering
    \includegraphics[width=0.5\linewidth]{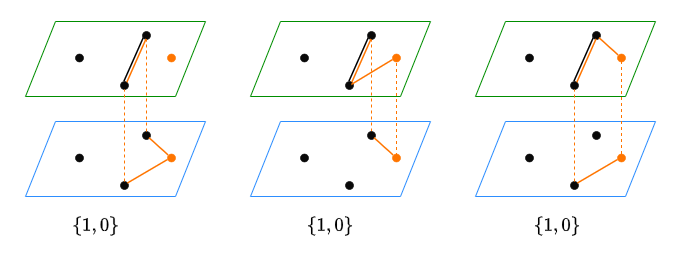}
    \caption{Dependent triangles in Case 4}
    \label{fig:32f1g1}
\end{figure}

\textbf{Case 5:} No new layers, two new nodes\\
Since we have two new nodes, only down edges can be shared. Once we have picked two layers out of $i,j,k$, the shared down edge is uniquely fixed and gives rise to four possible triangles $\beta$, see Figure~\ref{fig:32f0g2}. We have, again combining terms:
\[
B_5:= 3\binom{n}{3}\binom{n-3}{2} \sum_{(i,j,k)\in[L]^{3,\neq}} 4p_ip_j^2p_k^3q^4(1-q).
\]

\begin{figure}[H]
    \centering
    \includegraphics[width=0.66\linewidth]{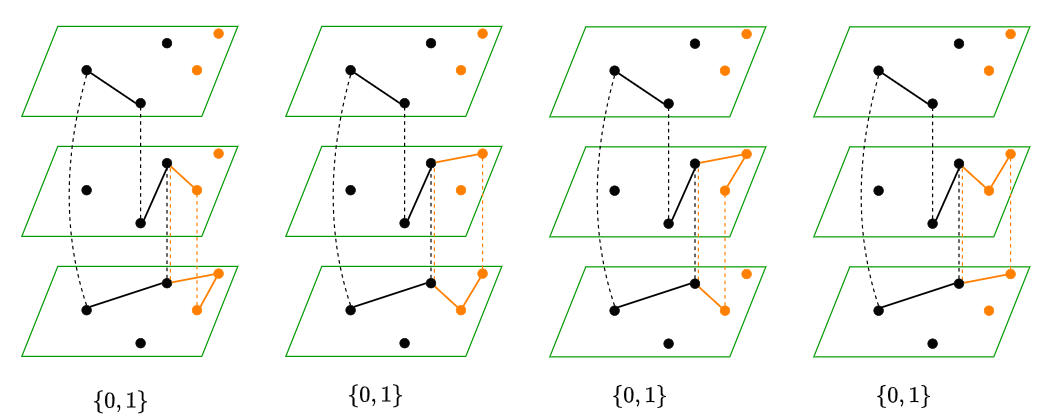}
    \caption{Dependent triangles for no new layers, two new nodes}
    \label{fig:32f0g2}
\end{figure}

Summing yields
\begin{align*}
    R_{3,2} \leq B_1+B_2 + B_3 +B_4 + B_5
    \leq& n^3 \sum_{(i,j,k)\in [L]^{3,\neq}} (2p_ip_j^2p_k q^4 + p_i^2p_j^3p_kq^4(1-q)) \\
    &+ \frac{3}{2}n^3 \sum_{(i,j,k,\ell)\in[L]^{4,\neq}} p_ip_jp_kp_\ell^2q^5 \\
    &+ n^4 \sum_{(i,j,k) \in [L]^{3,\neq}} (2p_i^2p_j^2p_kq^\red{4} + 2p_i^2p_j^3p_kq^4(1-q) + p_i^3p_jp_kq^4)\\
    &+ \frac{3}{2}n^4 \sum_{(i,j,k,\ell)\in[L]^{4,\neq}} p_ip_jp_k^2p_\ell q^5 \\
    &+ \frac{1}{4}n^5 \sum_{(i,j,k)\in[L]^{3,\neq}} 4p_ip_j^2p_k^3q^4(1-q).
\end{align*}

\subsubsection*{Bound for $R_{3,3}$}

Fix a 3D triangle $\alpha \in \Gamma_3$, and call the associated layers $i,j,k$. Here again we split the calculation into several subcases, depending on how many new nodes and layers we add.

\medskip \textbf{Case 1:} No new layers, no new nodes\\
Given three nodes and three layers, Figure~\ref{fig:33f0g0} shows all possible 3D triangles (in orange) we can create. Three of these configurations are of type $\{1,1\}$, in which case we have $\operatorname{Cov}(X_\alpha,X_\beta) \leq p_ip_j^2p_k^2q^5$ (if the shared edge is in layer $i$). This case contributes
\[
A_1 := 3\binom{n}{3} \sum_{i=1}^L\sum_{\substack{j=1\\j \neq i}}^L\sum_{\substack{k=1\\k \neq i,j}}^L p_ip_j^2p_k^2q^5 
\leq \frac{1}{2}n^3 \sum_{(i,j,k)\in [L]^{3,\neq}} p_ip_j^2p_k^2q^5
\]
to the bound on $R_{3,3}$.

\begin{figure}[H]
    \centering
    \includegraphics[width=1\linewidth]{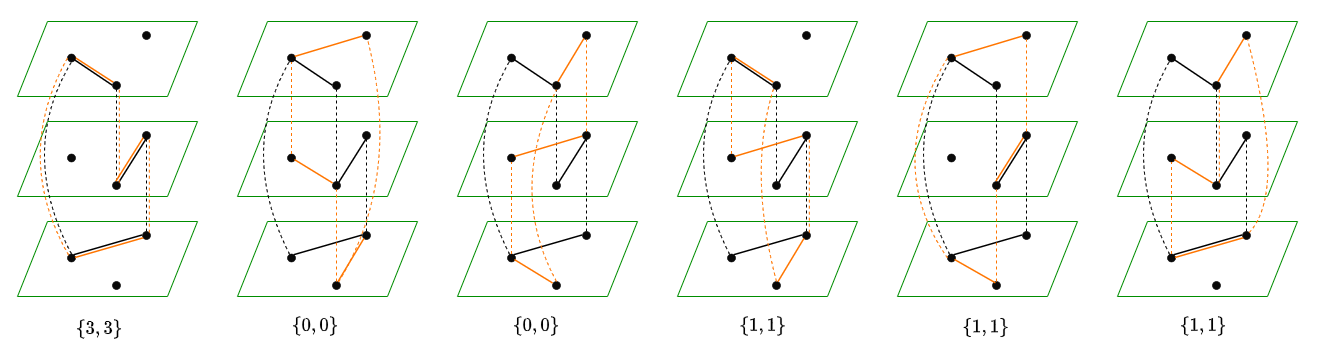}
    \caption{Dependent triangles with no new layer and no new nodes}
    \label{fig:33f0g0}
\end{figure}

\textbf{Case 2:} One new layer, no new nodes\\

To construct $\beta$, we can pick any two layers out of $i,j,k$, which gives $3$ options. Once this choice is fixed, Figure~\ref{fig:33f1g0} gives all possible triangles $\beta$ in orange. Here we add the following to the bound, grouping by symmetry as before:
\begin{align*}
    A_2:=& 3\binom{n}{3} \sum_{i=1}^L\sum_{\substack{j=1\\j \neq i}}^L\sum_{\substack{k=1\\k \neq i,j}}^L \sum_{\substack{\ell=1 \\ \ell \neq i,j,k}}^L (p_ip_jp_kp_\ell q^5 + 2p_i^2p_jp_kp_\ell q^6 + p_i^2p_j^2p_kp_\ell q^5(1-q))\\
    &\leq 3\binom{n}{3} \sum_{(i,j,k,\ell) \in [L]^{4,\neq}} (3p_ip_jp_kp_\ell q^5 + p_i^2p_j^2p_kp_\ell q^5(1-q))
\end{align*}
\begin{figure}[H]
    \centering
    \includegraphics[width=1\linewidth]{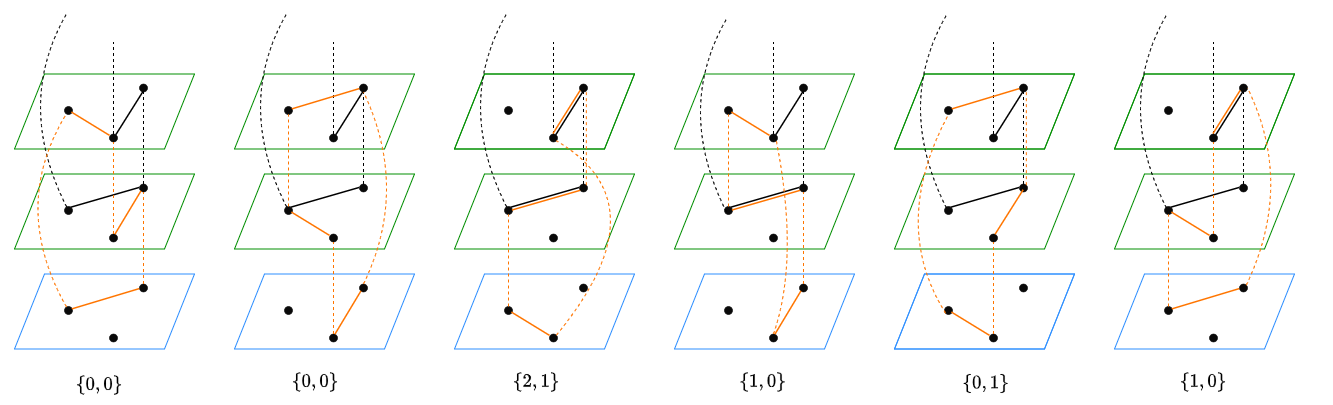}
    \caption{Dependent triangles with one new layer and no new nodes. Note that we omitted the top layer from the picture.}
    \label{fig:33f1g0}
\end{figure}

\textbf{Case 3:} No new layers, one new node\\
Given a new node (for which there are $n-3$ choices), there are $3$ choices which other two nodes the new one is connected to. Once this choice has been made, six triangles $\beta$ can be formed, see Figure~\ref{fig:33f0g1}. In this case we have
\begin{align*}
A_3 :&= 3\binom{n}{3}(n-3) \sum_{i=1}^L\sum_{\substack{j=1\\j \neq i}}^L\sum_{\substack{k=1\\k \neq i,j}}^L (2p_i^2p_j^2p_k^2 q^5(1-q) + p_i^2p_jp_k^2 q^4 + p_i^2p_jp_k^2 q^6)\\
&\leq 3\binom{n}{3}(n-3) \sum_{(i,j,k) \in [L]^{3,\neq}} (2p_i^2p_j^2p_k^2 q^5(1-q) + 2p_i^2p_jp_k^2 q^4).
\end{align*}
\begin{figure}[H]
    \centering
    \includegraphics[width=\linewidth]{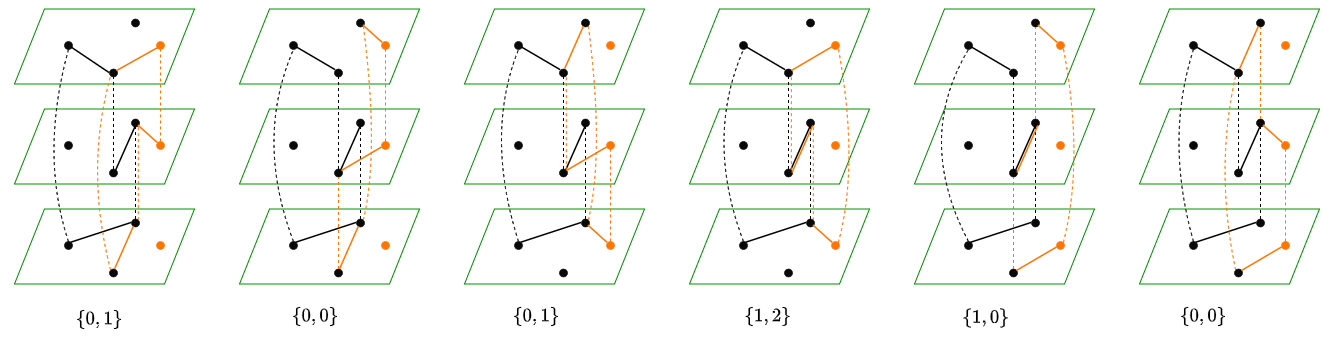}
    \caption{Dependent triangles for no new layers, one new node}
    \label{fig:33f0g1}
\end{figure}

\textbf{Case 4:} One new layer, one new node\\
Given a new layer $\ell$, there are $3$ ways to pick the other two layers among $i,j,k$. Once the choice of layers has been made, we need to pick two additional nodes among the three nodes of $\alpha$. \red{Two} of these choices results in shared edges between $\alpha$ and $\beta$, and this is illustrated in Figure~\ref{fig:33f1g1}. We thus add to the bound:
\begin{align*}
A_4:&= 3\binom{n}{3}(n-3) \sum_{i=1}^L\sum_{\substack{j=1\\j \neq i}}^L\sum_{\substack{k=1\\k \neq i,j}}^L \sum_{\substack{\ell=1 \\ \ell \neq i,j,k}}^L (\red{2}p_i^2p_j^2p_kp_\ell q^5(1-q) + \red{2}p_ip_j^2p_kp_\ell q^5 + \red{2}p_ip_j^2p_kp_\ell q^6) \\
&\leq 3\binom{n}{3}(n-3) \sum_{(i,j,k,\ell)\in [L]^{4,\neq}} (\red{2}p_i^2p_j^2p_kp_\ell q^5(1-q) + \red{4}p_ip_j^2p_kp_\ell q^5).
\end{align*}

\begin{figure}[H]
    \centering
    \includegraphics[width=1\linewidth]{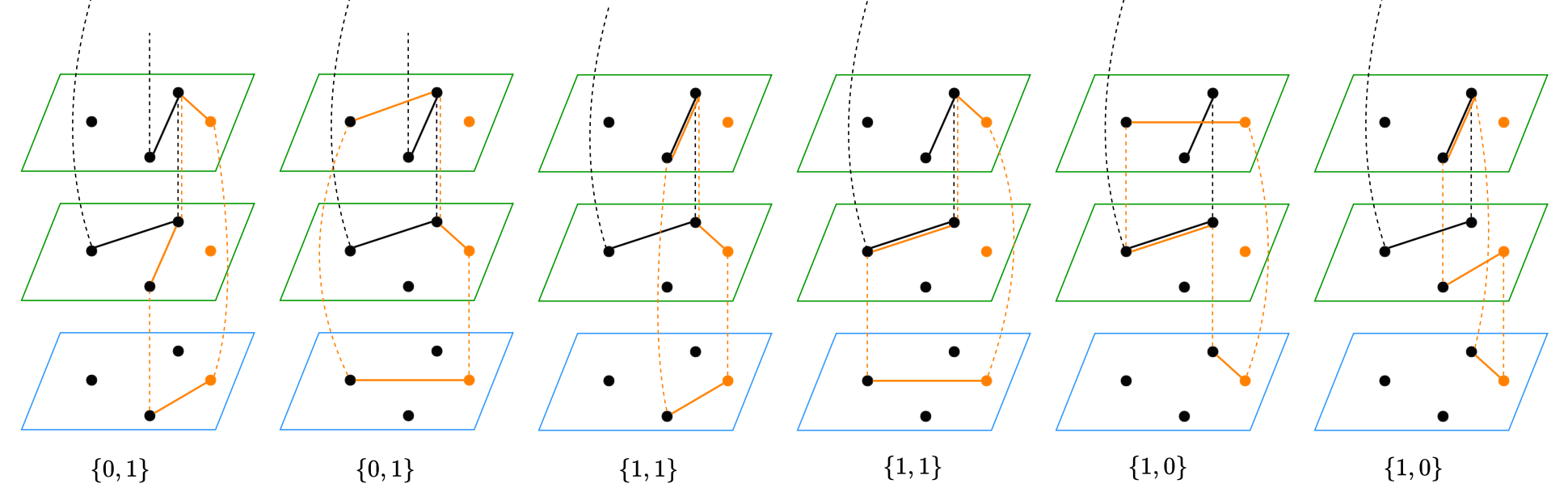}
    \caption{Dependent triangles for one new layer and one new node, omitting the top layer from the picture.}
    \label{fig:33f1g1}
\end{figure}

\textbf{Case 5:} Two new layers, no new nodes\\

Given two new layers $\ell,m$, we pick one of $i,j,k$ and the corresponding edge of $\alpha$ to get two possible triangles $\beta$, illustrated in Figure~\ref{fig:33f2g0}. Combining terms, we have thus
\[
A_5 := 3\binom{n}{3} \sum_{i=1}^L\sum_{\substack{j=1\\j \neq i}}^L\sum_{\substack{k=1\\k \neq i,j}}^L \sum_{\substack{\ell=1 \\ \ell \neq i,j,k}}^L \sum_{\substack{m=1 \\ m \neq i,j,k,\ell}}^L 2p_ip_jp_kp_\ell p_m q^6 \leq n^3 \sum_{(i,j,k,\ell,m) \in [L]^{5,\neq}} p_ip_jp_kp_\ell p_m q^6
\]

\begin{figure}[H]
    \centering
    \includegraphics[width=0.35\linewidth]{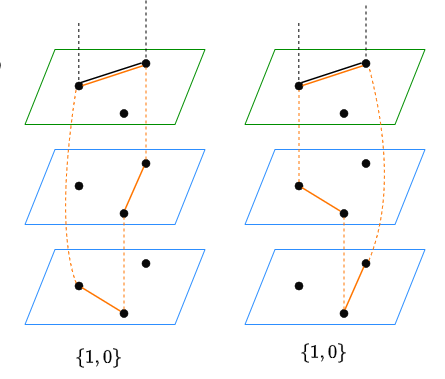}
    \caption{Dependent triangles for two new layers, no new nodes, omitting the top two layers.}
    \label{fig:33f2g0}
\end{figure}

\textbf{Case 6:} No new layers, two new nodes\\
Given two new nodes, we pick one of the three nodes of $\alpha$ and can form two triangles sharing a down edge, see Figure~\ref{fig:33f0g2}. We get
\[
A_6 := 3\binom{n}{3}\binom{n-3}{2} \sum_{i=1}^L\sum_{\substack{j=1\\j \neq i}}^L\sum_{\substack{k=1\\k \neq i,j}}^L 2p_i^2p_j^2p_k^2 q^5(1-q) \leq \frac{1}{2} n^5 \sum_{(i,j,k)\in[L]^{3,\neq}} p_i^2p_j^2p_k^2 q^5(1-q).
\]

\begin{figure}[H]
    \centering
    \includegraphics[width=0.35\linewidth]{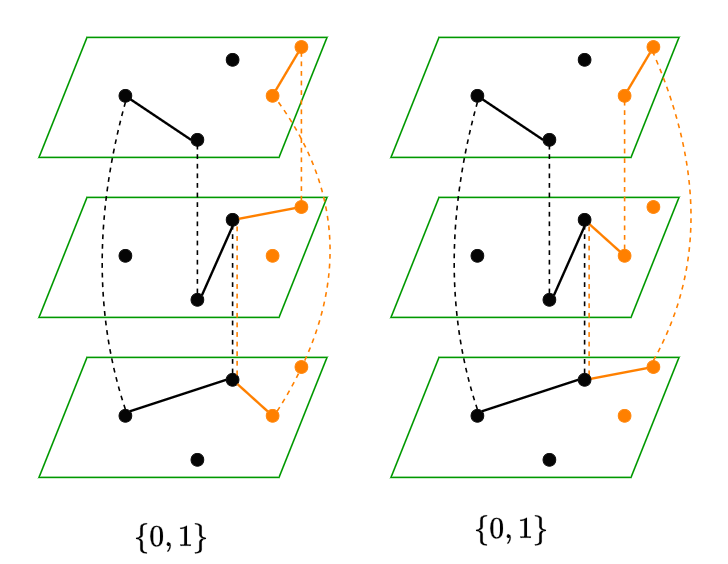}
    \caption{Dependent triangles for no new layers, two new nodes}
    \label{fig:33f0g2}
\end{figure}

\textbf{Case 7:} Two new layers, one new node\\
Given two new layers $\ell,m$, we select one of the layers $i,j,k$ to go with them. To create a dependent triangle $\beta$, we must choose the edge in that layer, which gives two possibilities for $\beta$, see Figure~\ref{fig:33f2g1}. We have
\[
A_7 := 3 \binom{n}{3} (n-3) \sum_{i=1}^L\sum_{\substack{j=1\\j \neq i}}^L\sum_{\substack{k=1\\k \neq i,j}}^L \sum_{\substack{\ell=1 \\ \ell \neq i,j,k}}^L \sum_{\substack{m=1 \\ m \neq i,j,k,\ell}}^L 2p_ip_jp_kp_\ell p_m q^6 \leq n^4 \sum_{(i,j,k,\ell,m) \in [L]^{5,\neq}} p_ip_jp_kp_\ell p_m q^6.
\]

\begin{figure}[H]
    \centering
    \includegraphics[width=0.35\linewidth]{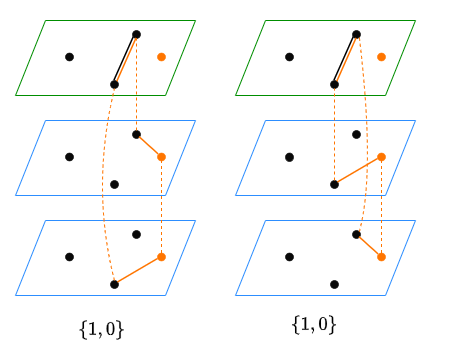}
    \caption{Dependent triangles for two new layers, one new node}
    \label{fig:33f2g1}
\end{figure}

\textbf{Case 8:} One new layer, two new nodes\\
Given one new layer $\ell$, we pick two out of the layers $i,j,k$ to form triangle $\beta$. Only down edges can be shared in this context because we have two new nodes, and once the choice of layers is fixed, the unique given down edge determines the choice of the third node. We have thus:
\[
A_8 := 3\binom{n}{3}\binom{n-3}{2} \sum_{(i,j,k,\ell)\in[L]^{4,\neq}} 2p_ip_j^2p_k^2p_\ell q^5(1-q).
\]
\begin{figure}[H]
    \centering
    \includegraphics[width=0.3\linewidth]{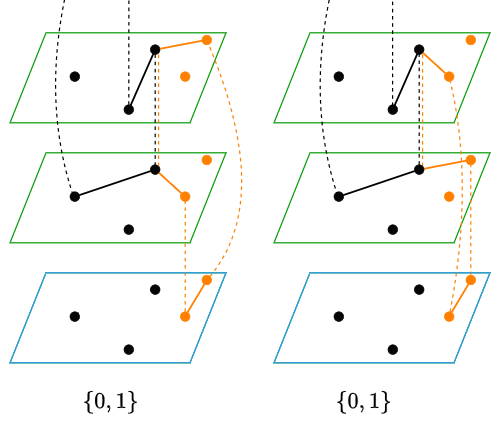}
    \caption{Dependent triangles for one new layer, two new nodes}
    \label{fig:33f1g2}
\end{figure}

Now, combining all of the above estimates, we get
\begin{align*}
R_{3,3} \leq A_1+\hdots+A_8
\leq& \frac{1}{2}n^3 \sum_{(i,j,k)\in [L]^{3,\neq}} p_ip_j^2p_k^2q^5 \\
&+ \frac{1}{2}n^3 \sum_{(i,j,k,\ell) \in [L]^{4,\neq}} (3p_ip_jp_kp_\ell q^5 + p_i^2p_j^2p_kp_\ell q^5(1-q)) \\
&+ \frac{1}{2}n^4 \sum_{(i,j,k) \in [L]^{3,\neq}} (2p_i^2p_j^2p_k^2 q^5(1-q) + 2p_i^2p_jp_k^2 q^4) \\
&+ \frac{1}{2}n^4 \sum_{(i,j,k,\ell)\in [L]^{4,\neq}} (\red{2}p_i^2p_j^2p_kp_\ell q^5(1-q) + \red{4}p_ip_j^2p_kp_\ell q^5) \\
&+ n^3 \sum_{(i,j,k,\ell,m) \in [L]^{5,\neq}} p_ip_jp_kp_\ell p_m q^6 \\
&+ \frac{1}{2} n^5 \sum_{(i,j,k)\in[L]^{3,\neq}} p_i^2p_j^2p_k^2 q^5(1-q) \\
&+ n^4 \sum_{(i,j,k,\ell,m) \in [L]^{5,\neq}} p_ip_jp_kp_\ell p_m q^6 \\
&+ \frac{1}{4} n^5 \sum_{(i,j,k,\ell)\in[L]^{4,\neq}} 2p_ip_j^2p_k^2p_\ell q^5(1-q). 
\end{align*}

This completes the proof.
\end{proof}

\begin{proof}[Proof of Theorem~\ref{theorem1}]

The proof of Theorem~\ref{theorem1} is now straightforward. We use the multivariate Poisson bounds from Theorem~\ref{thm:dtvmpoi} with 
the sets of indices $\Gamma=\Gamma_1 \cup \Gamma_2 \cup \Gamma_3$ introduced in Section~\ref{sec:background}. In our case,  for any index $\alpha \in \Gamma$, we have $\Gamma \setminus \{\alpha\} = \Gamma_\alpha^+$, since the presence of a triangle at $\alpha$ makes all triangles sharing an edge with a triangle at $\alpha$ more likely to appear, and does not influence any other triangles. Therefore $\Gamma_\alpha^-=\Gamma_\alpha^0=\emptyset$ and we have

\begin{align}
    \label{eq:tobound}
d_{T V}\left(\mathcal{L}\left(\left\{W_j\right\}_{j=1}^r\right),  \Pi_{j=1}^r \text{Po}(\lambda_j)\right) \leq \sum_{\alpha \in \Gamma} \mathbb{P}(X_\alpha=1)^2+\sum_{\alpha \in \Gamma} \sum_{\beta \in \Gamma\setminus\{\alpha\}} \operatorname{Cov}\left(X_\alpha, X_\beta\right).
\end{align}
 We have
\begin{align*}
    \sum_{\alpha \in \Gamma} \mathbb{P}(X_\alpha=1)^2 &= \sum_{\alpha = (\alpha_1^i,\alpha_2^i,\alpha_3^i) \in \Gamma_1} p_i^6 + \sum_{\alpha = (\alpha_1^i,\alpha_2^j,\alpha_3^j)\in \Gamma_2} p_i^2p_j^4q^\red{4} + \sum_{\alpha=(\alpha_1^i,\alpha_2^j,\alpha_3^k) \in \Gamma_3} p_i^\red{2}p_j^\red{2}p_k^\red{2}q^\red{6} \\
    &= \binom{n}{3} \sum_{i=1}^L p_i^6 + 3\binom{n}{3} \sum_{i=1}^L\sum_{\substack{j=1 \\ j \neq i}}^L p_i^2 p_j^4q^\red{4} + \binom{n}{3} \sum_{i=1}^L\sum_{\substack{j=1 \\ j \neq i}}^L \sum_{\substack{k=1 \\ k \neq i,j}}^L p_i^2p_j^2p_k^2q^6,
\end{align*}
which gives the first term in the assertion of Theorem~\ref{theorem1}. As for the second term, we have
\[
\sum_{\alpha \in \Gamma} \sum_{\beta \in \Gamma\setminus\{\alpha\}} \operatorname{Cov}\left(X_\alpha, X_\beta\right) = R_{1,1} + R_{2,2} + R_{3,3} + 2R_{2,1} + 2R_{3,1} + 2R_{3,2},
\]
where 
$
R_{i,j}= \sum_{\alpha \in \Gamma_i}\sum_{\substack{\beta\in \Gamma_j \\ \beta \neq \alpha}} \operatorname{Cov}(X_\alpha,X_\beta)$. 
Inserting this in \eqref{eq:tobound} yields the first assertion. 

\medskip

\red{The bound for the simpler case $p_i=p$ for all $i \in \{1,...,L\}$ and $q=1$, can be derived directly. Indeed, one has that
\[
\sum_{\alpha \in \Gamma} \mathbb{P}(X_\alpha=1)^2 = \binom{n}{3} L^3 p^6,
\]
since we first pick three nodes, then the three layers for the respective edges, each occurring with probability $p$. As for the term $\sum_{\alpha \in \Gamma}\sum_{\beta \in \Gamma_\alpha^+}\operatorname{Cov}(X_\alpha,X_\beta)$, note that there are $\binom{n}{3}L^3$ ways to pick the triangle $\alpha$. Once $\alpha$ has been fixed, there are $3$ possibilities for $\beta$ to share two edges with $\alpha$, with $L-1$ possible layers for the third edge of $\beta$. If $\beta$ shares only one edge and all nodes with $\alpha$, there are $3$ choices for the shared edge and $(L-1)^2$ choices for the layers of the other two edges. If $\beta$ shares exactly one edge and only two nodes with $\alpha$, then there are $n-3$ choices for the third node and $L^2$ possible layers for the other two edges. In total, we thus achieve
\begin{multline*}
\sum_{\alpha \in \Gamma}\sum_{\beta \in \Gamma_\alpha^+}\operatorname{Cov}(X_\alpha,X_\beta) \\
= L^3 \binom{n}{3} \big( 3(L-1)(p^4-p^6) + 3 (L-1)^2 (p^5-p^6) + 3(n-3)L^2(p^5-p^6)\big).    
\end{multline*}
Simplifying and using $p\leq 1$ and $L,n\geq 1$, we obtain
\[
d_{T V}\left(\mathcal{L}\left(W_1,W_2,W_3\right),  \prod_{j=1}^3 \text{Po}(\lambda_j)\right) \leq \frac{2}{3}n^3L^4p^4 + n^4L^5p^5,
\]
which concludes the proof of \eqref{eq:second}.
}
\end{proof}
\end{document}